  \def\obs#1{{\bf (*** #1 ***)} }

\def\obs#1{}     
\documentclass[twoside,letterpaper,draft,11pt]{amsart}
\usepackage{amsmath,amsthm,latexsym,xspace,amscd,amssymb,mathtools,color, enumerate, amsthm, amscd, amssymb, color,amsfonts}
\usepackage{mathrsfs}
\usepackage{relsize}

\DeclareUnicodeCharacter{00B4}{}

\usepackage[utf8]{inputenc}

\usepackage[mathscr]{eucal}
\usepackage{amssymb}
\usepackage[all]{xy}
\makeatletter\renewcommand\theenumi{\@roman\c@enumi}\makeatother

\newtheorem{teo1}{Theorem}[section]
\newtheorem{def1}[teo1]{Definition}
\newtheorem{lem1}[teo1]{Lemma}
\newtheorem{cor1}[teo1]{Corollary}
\newtheorem{prop1}[teo1]{Proposition}

\newtheorem{exe}[teo1]{Example}
\newtheorem{remark}[teo1]{Remark}

\newcommand{\m}{{}^{-1}}

\newcommand{\mcI}{{\mathcal{I}}}

\newcommand{\id}{{\rm id}}

\newcommand{\N}{{\mathbb N}}

\newcommand{\I}{{I}}

\newcommand{\A}{{A}}

\newcommand{\ot}{\otimes}
\newcommand{\af}{\alpha}
\newcommand{\bt}{\beta}

\newcommand{\tG}{\mathtt{G}}
\newcommand{\mcD}{\mathcal{D}}

\newcommand{\bI}{\mathbb{I}}

\newcommand{\G}{\mathcal{G}}

\def\ndv{\ {\mid \kern -0.7 em {\scriptstyle \not}} \ \ }

\def\nd{\ {\mid \kern -0.4 em {\scriptstyle \not}} \ \ }

\begin{document}

\thispagestyle{empty}

\title[Generalized matrix algebra]{Partial actions of groups on generalized matrix rings}

\author[D. Bagio]{Dirceu Bagio}
\address{ Departamento de Matem\'atica, Universidade Federal de Santa Catarina, 88040-900\\
	Florian\'opolis, Brasil}
\email{d.bagio@ufsc.br}


\author[H. Pinedo]{H\'ector Pinedo}
\address{Escuela de Matematicas, Universidad Industrial de Santander, Cra. 27 Calle 9 ´ UIS
	Edificio 45, Bucaramanga, Colombia}
\email{hpinedot@uis.edu.co}


\date{\today}
\begin{abstract} Let $n$ be a positive integer and $R=(M_{ij})_{1\leq i,j\leq n}$ be a generalized matrix ring. For each $1\leq i,j\leq n$, let $I_i$ be an ideal of the ring $R_i:=M_{ii}$ and denote $I_{ij}=I_iM_{ij}+M_{ij}I_j$. We give sufficient conditions for the subset $I=(I_{ij})_{1\leq i,j\leq n}$ of $R$ to be an ideal of $R$. Also, suppose that $\alpha^{(i)}$ is a partial action of a group $\mathtt{G}$ on $R_i$, for all $1\leq i\leq n$. We construct, under certain conditions, a partial action $\gamma$ of $\mathtt{G}$ on $R$ such that $\gamma$ restricted to $R_i$ coincides with $\alpha^{(i)}$. We study the relation between this construction and the notion of Morita equivalent partial group action  given in \cite{ADES}. Moreover, we investigate properties related to Galois theory for the extension $R^{\gamma}\subset R$. Some examples to illustrate the results are considered in the last part of the paper.
\end{abstract}

\maketitle

\setcounter{tocdepth}{1}

\section{Introduction} 

Let $n$ be a positive integer, $\{R_i\}_{i=1}^n$ a set of rings, $M_{ij}$ an  $(R_i,R_j)$-bimodule such that $M_{ii}=R_i$ and  $\theta_{ijk}:M_{ij}\otimes_{R_j} M_{jk}\to M_{ik}$ a bimodule morphism, for all  $1\leq i,j,k\leq n$.
An {\it $n$-by-$n$   generalized matrix ring $R$} is a square array
\begin{align*}
	R = \left( \begin{matrix} R_1 & M_{12} & \ldots & M_{1n} \\ 
		M_{21} & R_2 & \ddots & \vdots  \\ 
		\vdots & \ddots & \ddots & M_{(n-1)n} \\ 
		M_{n1} &\ldots & M_{n(n-1)} & R_n \end{matrix} \right) 
\end{align*}                         
Denote $m_{ij}m_{jk}:=\theta_{ijk}(m_{ij}\otimes m_{jk})\in M_{ik}$, for all $m_{ij}\in M_{ij}$ and $m_{jk}\in M_{kj}$. Assume that the bimodule morphisms $\theta_{ijk}$  satisfy the associativity relation
\begin{align*}
	(m_{ij}m_{jk})m_{kl}=m_{ij}(m_{jk}m_{kl}), \qquad i,j,k,l\in \bI_n:=\{1,\ldots,n\},
\end{align*}
and consider the addition on $R$ componentwise and the multiplication on $R$ row-column matrix multiplication.  In what follows, a generalized matrix ring $R$ as above will be denoted by $R=(M_{ij})_{i,j\in \bI_n}$. \vspace{.1cm}

Generalized matrix rings are a  large class of rings that  arise in 
several  areas of the theory of rings and modules  and have been extensively studied.  For instance, it is well known that   the class of $2$-by-$2$ generalized matrix rings  is in one-to-one correspondence with the Morita contexts; for details see \cite[Ch. 6]{F}.  Moreover, an arbitrary ring $A$ with a complete set of orthogonal idempotents $\{e_i\}_{i=1}^n$ has a Peirce decomposition, which can be arranged into an $n$-by-$n$ generalized matrix ring $A^{\pi}$ which is isomorphic to $A$.\vspace{.1cm}

Suppose that  $R=(M_{ij})_{i,j\in \bI_n}$ is a  $2$-by-$2$ generalized matrix ring. In  \cite{budanov}, the relationships between the lattices of ideals of the generalized matrix ring $R$  and the corresponding lattices of ideals of the rings $R_1$ and $R_2$ are described. Necessary and sufficient conditions for a pair of ideals $I$, $J$
of $R_1$ and $R_2$, respectively, to be the main diagonal of some ideal of the ring $R$ are also obtained in \cite{budanov}. More recently, multiplicative maps between generalized matrix rings were treated in \cite{FJ}.\vspace{.1cm}

On the other hand,  the notion of partial  group action appears as a weakening of the concept of global group action,  and has been 
applied in various areas of mathematics. An overview of the progress of partial action research can be seen in \cite{D}.  In the purely algebraic context, partial group actions appeared in \cite{DE}. Among other things, partial  group actions  were used to construct partial skew group rings (which is a generalization of skew group rings) as an algebraic analog of partial crossed product of $C^*$-algebras. The construction of partial skew group rings allowed the description of several important classes of rings and the introduction of the notion of partial Galois extension \cite{DFP}.\vspace{.1cm}


This work is structured as  follows. The background of partial action of groups, groupoids  and partial action of groupoids is presented in Section \ref{pl}. 
In Section \ref{gma} we work with generalized matrix rings. After recalling some basic facts about generalized matrix rings, we show that a certain class
of skew groupoid rings can be realized  as a generalized matrix rings. Also,  inspired by the results of \cite{budanov}, we give in Corollary \ref{cor:ideals} a sufficient condition for a finite family of ideals $I_j$ of $R_j$, $j\in \bI_n=\{1,\ldots,n\}$, to be the main diagonal of some ideal of the generalized matrix ring $R$. To do so, we introduce the notion of $\mathcal{I}$-symmetry for a generalized matrix ring $R$, where $\mathcal{I}=\{I_1,\ldots,I_n\}$. In order to construct a partial action of a group $\tG$ on the generalized matrix ring $R$, we define in Section 4 the notion of a datum $\mcD$ for $R$. The main result of this paper is Theorem \ref{teo-par}, in which we construct a partial action $\gamma=\gamma_{\mcD}$ of $\tG$ on $R$, for all datum $\mcD$ of $R$. Conditions for $\gamma$ to be unital, partial actions induced and Morita equivalent partial actions are explored in Section 4. In Section 5, we prove in Theorem \ref{teo-galois} that $R^{\gamma}\subset R$ is a partial Galois extension if and only if $R^{\alpha^{(i)}}_i\subset R_i$ is a partial Galois extension, for all $i\in \bI_n$. Section 6 is devoted to examples. In this last part, the main result is proved in Theorem \ref{teo-skew-iso}. Precisely, let 
$\af=(\A_g,\af_g)_{g\in \G}$ be a group-type partial action of a connected groupoid $\G$ on the ring $A=\oplus_{y\in \G_0}A_{y}$ and assume that $\G_0$ is finite. Then  the corresponding partial skew groupoid ring $ A*_{\af}\G$ is isomorphic to $R*_{\gamma} \tG$, where $R$ is a generalized matrix ring and $\gamma$ is a partial action of a group on $R$ obtained by the construction described in Section 4.

\subsection*{Notations and conventions}\label{subsec:conv}
By ring we mean an associative ring. The identity element of a group $\tG$ will be denoted by $e$ while $\id_{X}$ denotes the identity map of a set $X$. Let $I$ be an ideal of a unital ring $A$ and $M$ be a left $A$-module. We will denote by $IM$ the submodule of $M$ generated by the set $\{am\,:\,a\in I,\,m\in M\}$. Also, for a positive integer $n$, we will denote $\bI_n:=\{1,\ldots,n\}.$
The center of a ring $A$ will be denoted by $Z(A)$. For an $A$-module $M$, we will denote the action of $A$ on $M$ by concatenation, that is, $am$ means the action of $a\in A$ on $m\in M$.

\section{Preliminaries}\label{pl}

\subsection{Partial group actions}
We recall the definitions of partial group action  and the associated partial skew group ring. 

\begin{def1}\label{def-pa-group}{\rm Let $\tG$ be a group and $A$ be a ring. A {\it partial action} of $\tG$ on $A$ is a collection of pairs $\alpha=(\alpha_g,A_g)_{g\in \tG}$, where $A_g$ is an ideal of $A$ and $\alpha_g:A_{g\m}\to A_g$ is a ring isomorphism, for all $g\in \tG$,  such that:
		\begin{enumerate}[\rm (i)]
			\item $A_e=A$ and $\alpha_e=\id_{A}$,
			\item $\alpha_{h\m}(A_h\cap A_{g\m})\subset A_{(gh)\m}$, for all $g,h\in \tG$,
			\item $\alpha_g(\alpha_h(a))=\alpha_{gh}(a)$, for any $a\in \alpha\m_h(A_h\cap A_{g\m})$.
	\end{enumerate}}
\end{def1}

Other useful notions for our purposes in this paper are the following.

\begin{def1}\label{def-specific}
	{\rm A partial action $\af=(\A_g,\af_g)_{g\in \tG}$ of a group $\tG$ on a ring $A$ is called
		\begin{enumerate}[\rm (i)] 
			\item  {\it global}, if $A_g=A$, for all  $g\in \tG$,
			\item  {\it unital}, if each $A_g$ is a unital ring, that is, there exists a central idempotent element $1_g$ in $A$ such that $A_g=A1_g$, for all $g\in \tG$.
	\end{enumerate}} 
\end{def1}

Examples of partial actions of groups can be constructed as restrictions of global actions; see, for instance, \cite[Ex. 1.1]{DFP}. Conversely, under suitable assumptions, a partial action $\alpha$ is globalizable, that is, $\alpha$ is indeed a restriction of a global action. 
Precisely, by \cite[Thm. 4.5]{DE} a partial action $\alpha$ on a unital ring is globalizable if and only if $\alpha$ is unital. \smallbreak

Now we recall the notion of product partial action given in \cite[Def. 2.4]{ADES}.

\begin{def1}\label{def-prod-pa-group}{\rm Let $\tG$ be a group and $A$ be a ring. A {\it product partial action} of $\tG$ on $A$ is a collection of pairs $\alpha=(\alpha_g,A_g)_{g\in \tG}$, where $A_g$ is an ideal of $A$ and $\alpha_g:A_{g\m}\to A_g$ is a ring isomorphism, for all $g\in \tG$, that satisfies:
		\begin{enumerate}[\rm (i)]
			\item $A_e=A$ and $\alpha_e=\id_{A}$,
			\item $A_g^2=A_g$, $A_gA_h=A_h A_g$, for all $g,h\in \tG$,
			\item $\alpha_{g\m}(A_gA_{h})= A_{g} A_{gh}$, for all $g,h\in \tG$,
			\item $\alpha_g(\alpha_h(a))=\alpha_{gh}(a)$, for all $a\in A_{h\m}A_{(gh)\m}$.
	\end{enumerate}}
\end{def1}

As mentioned in \cite{ADES}, the discrepancy between Definition \ref{def-pa-group}  and Definition \ref{def-prod-pa-group} can be removed with the following concept.

\begin{def1}\label{par-act-reg}{\rm (\cite[Def. 2.5]{ADES})
A partial action $\alpha=(\alpha_g,A_g)_{g\in \tG}$ of a group $\tG$ on a ring $A$ is called {\it regular} if \[D_{g_1}\cap D_{g_2}\cap\ldots\cap D_{g_n}=D_{g_1} D_{g_2}\ldots D_{g_n},\qquad g_i\in \tG,\,\, i\in \bI_n,\,\, n>0.\]}
\end{def1}

\begin{remark}\label{ppa}{\rm Let $\alpha=(\alpha_g,A_g)_{g\in \tG}$ be a partial action of a group $\tG$ on a ring $A$. If $\alpha$ is unital, then it is clear that $\alpha$ is regular.}
\end{remark}

Let $\alpha=(\alpha_g,A_g)_{g\in \tG}$ be  a partial action of a group $\tG$ on a ring $A$. The {\it partial skew group ring} $A*_{\alpha}\tG$ is the set of all finite formal sums $\sum_{g\in \tG} a_g\delta_g$, where $a_g\in A_g$. The ring structure of  $A*_{\alpha}\tG$ is the following: addition is componentwise and multiplication is given by the formula
$$(a_g\delta_g)(a_h\delta_h)=\alpha_g(\alpha_{g\m}(a_g)a_h)\delta_{gh},\quad a_g\in A_g,\,a_h\in A_h,\,g,h\in \tG.$$
It follows from (ii) above that the multiplication on $A*_{\alpha}\tG$ is well-defined. In general, $A*_{\alpha}\tG$ is not associative. However, if the ideals $A_g$ of $A$ are generated by central idempotents $1_g$ (that is, $A_g=A1_{g}$) then $A*_{\alpha}\tG$ is associative.  We suggest \cite{DE} for more details on the associativity of $A*_{\alpha}\tG$. \smallbreak

We recall other definitions that will be useful throughout this work.

\begin{def1} {\rm Let $\alpha=(A, A_g)_{g\in \tG}$ and $\theta=(B, B_g)_{g\in \tG}$ be partial actions of a group $\tG$ on the rings $A$ and $B$, respectively. We say that:
		\begin{enumerate}[\rm (i)]
			\item\emph{ $\alpha$ and $\theta$ are equivalent}, if there is a ring isomorphism $\varphi: A\to B$ such that $\varphi (A_g)=B_g$,  for any $g\in G$, and $\varphi (\alpha_g(x))=\theta_g(\varphi(x))$, for every $ x\in A_{g\m}$;\vspace{.1cm}
			\item\emph {$\theta$ is an extension of $\alpha$}, if there exists a monomorphism $\iota:A\to B$ such that $\iota(A_g)\subset B_g$ and $\theta_g(\iota(x))=\iota(\alpha_g(x))$, for all $g\in \tG$ and $x\in A_{g\m}$.
	\end{enumerate}}
\end{def1}
\begin{remark}{\rm  It is straightforward to see that equivalent partial actions have isomorphic partial skew group rings.}
\end{remark}

\subsection{Partial groupoid actions}  We recall from  \cite{BP}   the notion of partial groupoid action on a ring as  a natural generalization of group partial action. For the convenience of the reader we start by giving some basic facts of groupoids.

A {\it groupoid} $\G$ is a small category in which every morphism is an isomorphism. The set of  objects of $\G$  is denoted by $\G_0$.
If $g:x\to y$ is a morphism of $\G$ then $s(g)=x$ and $t(g)=y$ are called the {\it source} and the {\it target} of $g$ respectively. Any object $x$ of $\G$ is identified with its identity morphism, that is, $x=\id_x$. The set of morphisms from $y$ to $z$ will be denoted by $\G(y,z)$ and $\G(x):=\G(x,x)=\{g\in \G:\,s(g)=t(g)=x\}$ is the {\it isotropy group} associated to an object $x$ of $\G$.\vspace{.05cm}

We will use the concatenation to denote the composition of morphisms of a groupoid $\G$. Hence, for $g,h\in \G$, there exists $gh$ if and only if $t(h)=s(g)$.
For each $g\in \G$, we have that $s(g)=g^{-1}g$ and $t(g)=gg\m$. Also, $s(gh)=s(h)$ and $t(gh)=t(g)$ for all $g,h\in \G$ with $t(h)=s(g)$.\vspace{.05cm}
The groupoid $\G$ is {\it connected} if $\G(y,z)\neq \emptyset$, for all $y,z\in \G_0$. Any groupoid is a disjoint union of connected subgroupoids. Indeed, an equivalence relation on $\G_0$ given by $x\sim y$ $\Leftrightarrow$ $\G(x,y)\neq \emptyset$  induces the decomposition 
\begin{equation}\label{direct-sum}
	\G=\dot\cup_{X\in \G_0/\!\sim}\G_X.
\end{equation}
Here, $\G_X$ is the full connected subgroupoid of $\G$ such that  $(\G_X)_0=X$. By construction, $\G$ is the disjoint union of the subgroupoids $\G_X$.

Let $X$ be a nonempty set and $X^2=X\times X$. Then $X^2$ is a groupoid, called the {\it coarse groupoid} of $X.$ The source and target maps of $X^2$ are, respectively, $s(x,y)=x$ and $t(x,y)=y$, for all $x,y\in X$. The rule of composition is given by: 
\begin{equation}\label{coarse}(y,z)(x,y)=(x,z),\text{    for all      } x,y,z\in X.\end{equation}  It is well known that the direct product $\G\times\mathcal{H}$ of groupoids $\G$ and $\mathcal{H}$ is a groupoid with source and target maps given by
\begin{equation}\label{st} s(g,h)=(s(g),s(h))\,\, \text{    and    }\,\,t(g,h)=(t(g),t(h)),\end{equation} for all $g\in \G$ and $h\in \mathcal{H}$. Given composable morphisms $(g,h),(g',h')\in \G\times \mathcal{H}$, the product is defined by $(g,h)(g',h')=(gg',hh')$. The next result about the structure of connected groupoids is well known and its proof will be omitted; see, for instance,  \cite[Prop 2.1]{BPP}. 

\begin{prop1}\label{group:connec}
Let $\G$ be a connected groupoid, $x\in \G_0$ and $h_y\in \G(x,y)$, for each $y\in \G_0$, with $h_x=x$. Then the map $\psi:\G \to  \G_0^2\times \G(x)$ given by
\[\psi(g)=\left((s(g),t(g)),g_x\right),\qquad g_x:=h\m_{t(g)}gh_{s(g)}\in \G(x),\]
is an isomorphism of groupoids.
\end{prop1}

\begin{def1}\label{def-pa-groupoid}
	{\rm Let $\G$ be a groupoid and $A$ be a ring. A {\it partial action} of $\G$ on $A$ is a collection of pairs $\alpha=(\alpha_g,A_g)_{g\in \tG}$, where  $\A_g$ is an ideal of $\A_{t(g)}$, $\A_{t(g)}$ is an ideal of $\A$ and $\af_g:\A_{g\m}\to \A_g$ is a ring isomorphism, for all $g\in \G$,  that satisfies:
\begin{enumerate}[\rm (i)]
	\item $\alpha_x=\id_{A_x}$, for all $x\in \G_0$,
	\item $\alpha\m_h(A_h\cap A_{g\m})\subset A_{(gh)\m}$, for all $g,h\in \G$ such that $t(h)=s(g)$,
	\item $\alpha_g(\alpha_h(a))=\alpha_{gh}(a)$, for all $a\in \alpha\m_h(A_h\cap A_{g\m})$ and $g,h\in \G$ such that $t(h)=s(g)$.
\end{enumerate}	}
\end{def1}

Notice that $\af_{(x)}=(\A_h,\af_h)_{h\in\G(x)}$, which is obtained by restriction of $\alpha$, is a partial action of the isotropy group $\G(x)$ on the ring $\A_x$, for each $x\in\G_0$.\smallbreak

Since we will work with some specific types of partial actions of groupoids we recall the following.

\begin{def1}\label{def-specific-groupoid}
	{\rm A partial action $\af=(\A_g,\af_g)_{g\in \G}$ of a groupoid $\G$ on a ring $A$ is called
		\begin{enumerate}[\quad\rm (a)] 
			\item  {\it global}, if $A_g=A_{t(g)}$, for all $g\in \G$,
			\item  {\it unital}, if each $A_g$ is a unital ring, that is, there exists a central idempotent element $1_g$ in $A$ such that $A_g=A1_g$, for all $g\in \G$,
			\item {\it group-type}, if there exist $x\in \G_0$ and a set $\{h_y\,:\,h_y\in \G(x,y)\}$ of morphisms of $\G$ such that $h_x=x$ and
			\begin{align*} 
				A_{h\m_y}=A_x \ \ \text{and} \ \ A_{h_y}=A_y, \ \ \text{ for all } \ y\in\G_0.
			\end{align*}
	\end{enumerate}} 
\end{def1}

It is clear that any global groupoid action is group-type. The notion of group-type partial action does not depend on the choice of object $x$, see \cite[ Rmk. 3.4]{BPP}.  
As in the group case, examples of groupoid partial actions are constructed from groupoid global actions. The globalization problem for partial groupoid actions  on unital rings was studied in \cite{BP} and it was proved that a partial groupoid action $\alpha$ on a unital ring  is globalizable if and only if $\alpha$ is unital. \smallbreak

Let $\alpha=(\alpha_g,A_g)_{g\in \G}$ be  a partial action of a group $\G$ on a ring $A$. The {\it partial skew groupoid ring} $A\star_\af\G$ associated to $\alpha$ is the set of finite formal sums $\sum_{g\in \G}a_g\delta_g$, where $a_g\in A_g$, with the usual addition and multiplication induced by the following rule
\[
(a_g\delta_g)(a_h\delta_h)=
\begin{cases}
	\alpha_g(\alpha_{g\m}(a_g)a_h)\delta_{gh}, &\text{if $s(g)=t(h),$}\\
	0,  &\text{otherwise},
\end{cases}
\]
for all $g, h\in \G$, $a_g\in A_g$ and $a_h\in A_h$.
If the ideals $A_g$ are unital, then the partial skew groupoid ring $A\star_\af\G$ is an associative ring. If in addition $\G_0$ is finite, then $A\star_\af\G$ is  unital with identity element $1_{A\star_\af\G}=\sum_{y\in\G_0}1_y\delta_y$, where $1_y$ is the central idempotent of $A$ such that $A_y=A1_{y}$ (see $\S\,3$ of \cite{BFP} for more details). When $\alpha$ is global we say that $A\star_\af\G$ is a skew groupoid ring. \smallbreak

Let $\af=(\A_g,\af_g)_{g\in \G}$ be a group-type partial action of a connected groupoid $\G$ on the ring $A=\oplus_{y\in \G_0}A_{y}$. Thus, there exist $x\in \G_0$ and a set $\{h_y\,:\,h_y\in \G(x,y)\}$ of morphisms of $\G$ such that $h_x=x$ and $\alpha_{h_y}:A_x\to A_y$, for all $y\in \G_0$.
We will assume that $\G_0=\{x=x_1,x_2,\ldots,x_n\}$ is finite. According to \cite[Lemma 4.1]{BPP}, we have the following global action of $\G_0^2=\G_0\times\G_0$ on $A$:
\begin{align}\label{action-beta}
	&\bt=(A_{(x_i,x_j)},\bt_{(x_i,x_j)})_{(x_i,x_j)\in \G_0^2}, &&A_{(x_i,x_j)}=A_{x_j},&& \beta_{(x_i,x_j)}=\af_{h_{x_j}}\af_{h\m_{x_i}}.&
\end{align}

Moreover, by  \cite[Lem. 4.3]{BPP}, the pair $\varepsilon=(C_g, \varepsilon_g)_{g\in \tG}$ is a partial action of the isotropy group $\tG=\G(x)$ on $C=A*_{\beta}\G_0^2$, where
\begin{align}\label{pa-oa}
C_g=\mathlarger{\mathlarger{\oplus}}_{u\in\G_0^2}\,\,\af_{h_{t(u)}}(A_g)\delta_u,\qquad \varepsilon_g(\af_{h_{t(u)}}(a)\delta_u)= \af_{h_{t(u)}}(\af_g(a))\delta_u,	
\end{align}
for all $a\in A_{g\m}$ and $u\in \G_0^2$. It follows from   \cite[Thm 4.4]{BPP} that the map
\begin{align}\label{teo-skew}
\psi: A*_{\af}\G\to (A*_{\beta}\G_0^2)*_{\varepsilon} \tG,\qquad a\delta_g\mapsto\big( a\delta_{(s(g),t(g))}\big)\delta_{h\m_{t(g)}gh_{s(g)}},
\end{align}
is a ring isomorphism. 

\section{Generalized matrix rings}\label{gma}

Let $n$ be a positive integer, $\{R_i\}_{i=1}^n$ a set of rings, $M_{ij}$ an  $(R_i,R_j)$-bimodule such that $M_{ii}=R_i$ and  $\theta_{ijk}:M_{ij}\otimes_{R_j} M_{jk}\to M_{ik}$ a bimodule morphism, for all  $1\leq i,j,k\leq n$.
We recall that an {\it $n$-by-$n$   generalized matrix ring $R$} is a square array
	\begin{align}\label{def-gmr}
R = \left( \begin{matrix} R_1 & M_{12} & \ldots & M_{1n} \\ 
	                         M_{21} & R_2 & \ddots & \vdots  \\ 
	                         \vdots & \ddots & \ddots & M_{(n-1)n} \\ 
	                         M_{n1} &\ldots & M_{n(n-1)} & R_n \end{matrix} \right).
\end{align}                         
For $m_{ij}\in M_{ij}$ and $m_{jk}\in M_{kj}$, we denote $m_{ij}m_{jk}:=\theta_{ijk}(m_{ij}\otimes m_{jk})\in M_{ik}$. The bimodule morphisms $\theta_{ijk}$ must satisfy the associativity relation
\begin{align}\label{rel-assoc}
(m_{ij}m_{jk})m_{kl}=m_{ij}(m_{jk}m_{kl}), \qquad i,j,k,l\in \bI_n.	
\end{align}
Notice that  $\theta_{iii}$ is determined by the ring multiplication in $R_i$, while $\theta_{ijj}$ and $\theta_{jjk}$ are determined by the bimodule structures. The addition on $R$ is  componentwise and multiplication on $R$ is row-column matrix multiplication. It is clear that if $R_i$ is a unital ring for all $i\in \bI_n$, then $R$ is a unital ring and the identity of $R$ is the diagonal matrix $1_R=\operatorname{diag}\big(1_{R_1},\ldots,1_{R_n}\big)$. In what follows, a generalized matrix ring $R$ as in \eqref{def-gmr} will be denoted by $R=(M_{ij})_{i,j\in \bI_n}$.

\begin{remark}\label{sm1} {\rm Let $R=(M_{ij})_{i,j\in \bI_n}$ be a generalized matrix ring as in \eqref{def-gmr} and assume that $R_i$ is a unital ring, for all $i\in \bI_n$. Notice that $\Upsilon_{ij}=(R_i, R_j, M_{ij}, M_{ji}, \theta_{iji}, \theta_{jij})$ is a Morita context, for all $i,j\in \bI_n$. If $\theta_{iji}$ and  $\theta_{jij}$ are surjective bimodule morphisms (that is, the previous Morita context is strict), then $ M_{ij}$ and  $M_{ji}$ are finitely generated projective modules (both right and left) \cite[pg.167]{Ja}. 
}
\end{remark} 

\begin{exe}{\rm Let $(A,B,M,N,\mu,\nu)$ be a Morita context. Then 
$$ R = \left( \begin{matrix} A & M \\ N & B \end{matrix} \right)$$
is a $2$-by-$2$  generalized matrix ring which is called a Morita ring.
Conversely, it is clear that each $2$-by-$2$  generalized matrix ring determines a Morita context. Also, if $R_i=S$, $M_{ij}=S$ and $\theta_{ijk}=\id_{S}$, for all $1\leq i,j,k\leq n$, then the $n$-by-$n$ generalized matrix ring $R$ defined in \eqref{def-gmr} coincides with the usual matrix ring $\operatorname{M}_n(S)$. }
\end{exe}

\subsection{Skew groupoid ring as generalized matrix ring}\label{subsec-skew}

Let $n$ be a positive integer and consider the groupoid $\G=\mathbb{I}_n\times \mathbb{I}_n$ with partial operation given by $(j,k)(i,j)=(i,k)$. Notice that $\G_0=\{(i,i)\,:\,i\in \mathbb{I}_n\}$ is  the set of objects of $\G$ and the source and target maps $s,t:\G\to \G_0$ are given respectively by $s(i,j)=(i,i)$ and $t(i,j)=(j,j)$.

Suppose that $\{R_i\}_{i=1}^n$ is a set of unital rings. Assume that there are ring isomorphisms $\theta_i:R_1\to R_i$, for all $i\in \mathbb{I}_n$. Then, we define $\theta_{ij}:R_i\to R_j$ by $\theta_{ij}=\theta_j \theta\m_i$.
Notice that $\theta_{ii}=\id_{R_i}$ and $\theta_{jk}\theta_{ij}=\theta_{ik}$, for all $i,j,k\in \mathbb{I}_n$. In this case, we can define an action $\theta$ of $\G$ on $A=\prod_{i\in \mathbb{I}_n} R_i$ in the following way. Take the ideals $A_{(i,j)}:=R_{j}$ and the ring isomorphisms given by $\theta_{(i,j)}=\theta_{ij}:R_i\to R_j$. It is clear that $\theta=(A_u,\theta_u)_{u\in \G}$ is an action of $\G$ on $A$. Then we can consider the skew groupoid ring $A*_{\theta}\G=\oplus_{i,j\in \mathbb{I}_n} A_{j}\delta_{(i,j)}$, where $\delta_{(i,j)}$ are symbols. Recall that the sum on $A*_{\theta} \G$ is componentwise and the multiplication is defined by: 
\[(a\delta_{(k,l)})(b\delta_{(i,j)})=\begin{cases} a\theta_{jl}(b)\delta_{(i,l)},& \text{if } j=k,\\
	0,& \text{otherwise.}
\end{cases}\]

\begin{prop1}\label{prop:skew-gmr} Assume $\G$, $A$ and $\theta$ as above. Then the skew groupoid ring $A*_{\theta} \G$ is isomorphic to the following $n$-by-$n$ generalized matrix ring
	\begin{align}\label{prop-gmr}
		R =\left( \begin{matrix} R_1\delta_{(1,1)}&  R_1\delta_{(2,1)} & \ldots &  R_1\delta_{(n,1)} \\ 
			R_2\delta_{(1,2)} & R_2\delta_{(2,2)} & \ldots &  R_2\delta_{(n,2)}  \\ 
			\vdots & \ddots & \ddots & \vdots \\ 
			R_n\delta_{(1,n)}&  R_n\delta_{(2,n)} &\ldots  & R_n\delta_{(n,n)} \end{matrix} \right).
	\end{align}
\end{prop1}

\begin{proof} Let $i,j\in \mathbb{I}_n.$
	Observe that $R_i\delta_{(i,i)}$ is a ring. By setting $M_{ij}=R_i\delta_{(j,i)}$, one would have  that $M_{ij}$ is an $(R_i\delta_{(i,i)},R_j\delta_{(j,j)})$-bimodule with structures
	\begin{align}\label{eq:structure-module-skew}
	r\delta_{(i,i)}\cdot s\delta_{(j,i)}=rs\delta_{(j,i)},\qquad  s\delta_{(j,i)}\cdot t\delta_{(j,j)}=s\theta\m_{ij}(t)\delta_{(i,j)},\qquad r,s\in R_i,\,\,t\in R_j.	
	\end{align}
	Also, consider the following bimodule maps: for all $r\in R_i$ and $s\in R_j$,
	\begin{align}\label{one-more}
	\theta_{ijk}:M_{ij}\otimes_{R_j} M_{jk}\to M_{ik},\qquad \theta_{ijk}\big(r\delta_{(j,i)}\otimes s\delta_{(k,j)}\big)=r\theta\m_{ij}(s)\delta_{(k,i)}.	
	\end{align}
	Since $\theta_{jk}\theta_{ij}=\theta_{ik}$, it follows that $(m_{ij}m_{jk})m_{kl}=m_{ij}(m_{jk}m_{kl})$, for all $m_{ij}\in M_{ij}$, $m_{jk}\in M_{jk}$ and $m_{kl}\in M_{kl}$.
	Thus $R=(M_{ij})_{i,j\in \bI_n}$ is a $n$-by-$n$ generalized matriz ring. 
	Now, consider the map $\varphi$ from $A*_{\theta} \G$ to $R$ given by 
	\[\varphi\Big(\sum_{i,j=1}^{n}a_{ij}\delta_{(j,i)}\Big)= \left( \begin{matrix} a_{11}\delta_{(1,1)}&  a_{12}\delta_{(2,1)} & \ldots &  a_{1n}\delta_{(n,1)} \\ 
		\vdots & \vdots & \vdots & \vdots \\ 
		a_{n1}\delta_{(1,n)}&  a_{n2}\delta_{(2,n)} &\ldots  & a_{nn}\delta_{(n,n)} \end{matrix} \right),\quad a_{ij}\in R_i. \]
	It is straightforward to check that $\varphi$ is a ring isomorphism.
\end{proof}

\subsection{Ideals of generalized matrix rings} Let $n$ be a positive integer and consider a generalized matrix ring $R=(M_{ij})_{i,j\in \bI_n}$ as in \eqref{def-gmr}. For each $j\in \bI_n$, we take an ideal $I_j$ of $R_j$ and we set
\begin{align}\label{form-for-ideals}
	I_{jk}:=I_jM_{jk}+M_{jk}I_k, \qquad  i,j,k\in \bI_n.
\end{align} 
Observe that $I_{jk} \subset M_{jk}$. Also, if $R_j$ is unital then $I_{jj}=I_j$.

\begin{def1}\label{def-symetric} Let $n$ be a positive integer, $R=(M_{ij})_{i,j\in \bI_n}$ be a generalized matrix ring as in \eqref{def-gmr}, and $\mcI=\{I_j\subset R_j\,:\,I_j \text{ is an ideal of }R_j,\text{ for all }j\in \bI_n\}$ be a finite family of ideals. The $(R_i,R_j)$-bimodule $M_{ij}$ will be called $\mcI$-symmetric if $I_iM_{ij}=M_{ij}I_j$. When $M_{ij}$ is $\mcI$-symmetric for all $i,j\in \bI_n$, we say that $R$ is $\mcI$-symmetric.
\end{def1}

\begin{exe}\label{ex-symmetric} {\rm Let $n$ be a positive integer, $\G=\mathbb{I}_n\times \mathbb{I}_n$ be the groupoid considered in $\S\,$\ref{subsec-skew} and $\{R_i\}_{i=1}^n$ be  a set of unital rings. Also, let $A=\prod_{i\in \mathbb{I}_n} R_i$ and $\theta=(A_u,\theta_u)_{u\in \G}$ be the partial action of $\G$ on $A$ defined in $\S\,$\ref{subsec-skew}. By Proposition \ref{prop:skew-gmr}, the partial skew groupoid ring $A\star_{\theta}\G$ is isomorphic to the generalized matrix ring $R=(M_{ij})_{i,j\in\bI_n}$, where $M_{ij}=R_i\delta_{(j,i)}$. Take an ideal $I_1$ of $R_1$. Then  $I_j:=\theta_j(I_1)$ is an ideal of $R_j$, for all $j\in \bI_n$. Hence, $I_j\delta_{(j,j)}$ is an ideal of $R_j\delta_{(j,j)}$ and we consider $\mcI=\{I_j\delta_{(j,j)}\,:\,j\in \bI_n\}$. It follows from \eqref{eq:structure-module-skew} that 
		\[I_j\delta_{(j,j)}R_j\delta_{(k,j)}=\theta_{j}(I_1)\delta_{(j,j)}R_j\delta_{(k,j)}=\theta_j(I_1)\delta_{(k,j)}\]
		and 
		\begin{align*}
			R_j\delta_{(k,j)}I_k\delta_{(k,k)}=R_j\delta_{(k,j)}\theta_k(I_1)\delta_{(k,k)}=R_j\theta\m_{jk}(\theta_k(I_1))\delta_{(k,j)}=\theta_j(I_1)\delta_{(k,j)}.
		\end{align*}
		Thus $I_j\delta_{(j,j)}M_{jk}=M_{jk}I_k\delta_{(k,k)}$, for all $j,k\in \bI_n$. Hence, $R$ is $\mcI$-symmetric.}
\end{exe}

The next result characterizes when a generalized matrix ring $R=(M_{ij})_{i,j\in \bI_n}$ is $\mcI$-symmetric for the case where $R_i=M_{ii}$ is unital, for all $i\in \bI_n$.

\begin{prop1}\label{prop-symmetric} Let $n$ be a positive integer, $R=(M_{ij})_{i,j\in \bI_n}$ be a generalized matrix ring as in \eqref{def-gmr}, and $\mcI=\{I_j\subset R_j\,:\,I_j \text{ is an ideal of }R_j,\text{ for all }j\in \bI_n\}$ be a finite family of ideals. Assume that $R_i=M_{ii}$ is unital, for all $i\in \bI_n$. The following statements are equivalent:
\begin{enumerate}[\rm (i)]
	\item $R$ is $\mcI$-symmetric,\vspace{.1cm}
	\item $M_{ij}I_jM_{jk}\subset M_{ik}I_k$ and $M_{ij}I_jM_{jk}\subset I_iM_{ik}$, for all $i,j,k\in \bI_n$.
\end{enumerate}	
\end{prop1}

\begin{proof}
Since $M_{ij}$ and $M_{jk}$ are $\mcI$-symmetric, it follows that
\[M_{ij}I_jM_{jk}=M_{ij}M_{jk}I_k\subset M_{ik}I_k,\qquad M_{ij}I_jM_{jk}=I_iM_{ij}M_{jk}\subset I_iM_{ik}.\]
Conversely, assume that (ii) is true. Since $R_i$ is unital, we have that
\[I_iM_{ij}=R_iI_iM_{ij}=M_{ii}I_iM_{ij}\subset M_{ij}I_j. \]
Also, $M_{ij}I_j=M_{ij}I_jR_j=M_{ij}I_jM_{jj}\subset I_iM_{ij}$. Thus $I_iM_{ij}=M_{ij}I_j$ and hence $M_{ij}$ is $\mcI$-symmetric.
\end{proof}

Let $n$ be a positive integer, $R=(M_{ij})_{i,j\in \bI_n}$ be a generalized matrix ring as in \eqref{def-gmr} and $\mcI=\{I_j\subset R_j\,:\,I_j \text{ is an ideal of }R_j,\text{ for all }j\in \bI_n\}$ a finite family of ideals. In order to construct an ideal of $R$, we consider the subset $I$ of $R$ given by
\begin{align}\label{expl-ideal}
	I = \left( \begin{matrix} I_{11} & I_{12} & \ldots & I_{1n}\\ 
		I_{21} & I_{22} & \ldots &  I_{2n}  \\ 
		\vdots & \vdots & \vdots &\vdots \\ 
		I_{n1} &I_{n2} &\ldots  & I_{nn}\end{matrix} \right),
\end{align} 
where $I_{jk}$ was defined in \eqref{form-for-ideals}. Then we have the following result.

\begin{prop1}\label{prop:ideals-gm} Let $n$ be a positive integer,  $R=(M_{ij})_{i,j\in \bI_n}$ be the generalized matrix ring given by \eqref{def-gmr}, $\mcI=\{I_j\subset R_j\,:\,I_j \text{ is an ideal of }R_j,\text{ for all }j\in \bI_n\}$ be a finite family of ideals, and $I$ be the subset of $R$ defined in \eqref{expl-ideal}. Then the following statements hold.
\begin{enumerate}[\rm (i)]
	\item If $M_{ij}I_jM_{jk}\subset M_{ik}I_k$, for all $i,j,k\in \I_n$, then $I$ is an ideal of $R$.\vspace{.15cm}
	\item If $M_{ij}I_jM_{jk}\subset I_iM_{ik}$, for all $i,j,k\in \I_n$, then $I$ is an ideal of $R$.
\end{enumerate}	
\end{prop1}

\begin{proof}
Let $i,j,k\in \bI_n$. Notice that  $I$ is an ideal of $R$ if and only if  $M_{ij}I_{jk}\subset I_{ik}$ and $I_{ij}M_{jk}\subset I_{ik}$. Suppose that $M_{ij}I_jM_{jk}\subset M_{ik}I_k$. Then
\[M_{ij}I_{jk}=M_{ij}I_jM_{jk}+M_{ij}M_{jk}I_k\subset M_{ik}I_k\subset I_{ik},\]
and
\[I_{ij}M_{jk}=I_{i}M_{ij}M_{jk}+M_{ij}I_jM_{jk}\subset I_iM_{ik}+M_{ik}I_k=I_{ik},\]
which implies (i). The proof of (ii) is similar.
\end{proof}

\begin{cor1} \label{cor:ideals}  Let $n$ be a positive integer,  $R=(M_{ij})_{i,j\in \bI_n}$ the generalized matrix ring given by \eqref{def-gmr}, $\mcI=\{I_j\subset R_j\,:\,I_j \text{ is an ideal of }R_j,\text{ for all }j\in \bI_n\}$ a finite family of ideals and $I$  the subset of $R$ defined in \eqref{expl-ideal}. If $R$ is $\mcI$-symmetric then $I$ is an ideal of $R$.
\end{cor1}	
\begin{proof}
Since $R$ is $\mcI$-symmetric it follows from Proposition \ref{prop-symmetric} that $M_{ij}I_jM_{jk}\subset M_{ik}I_k$ and $M_{ij}I_jM_{jk}\subset I_iM_{ik}$, for all $i,j,k\in \I_n$. then $I$ is an ideal of $R$. Then, by Proposition \ref{prop:ideals-gm}, the result follows.
\end{proof}

\begin{exe} \label{ex-ideal-skew}\rm{Let $n$ be a positive integer, $\G=\mathbb{I}_n\times \mathbb{I}_n$ be the groupoid considered in $\S\,$\ref{subsec-skew}, $\{R_i\}_{i=1}^n$ a set of unital rings,  $A=\prod_{i\in \mathbb{I}_n} R_i$ and $\theta=(A_u,\theta_u)_{u\in \G}$ the partial action of $\G$ on $A$ defined in $\S\,$\ref{subsec-skew}. By Proposition \ref{prop-gmr}, $A*_{\theta} \G$ is isomorphic to the generalized matrix ring $R=(M_{ij})_{i,j\in\bI_n}$, where $M_{ij}=R_i\delta_{(j,i)}$. Let $\mcI$ be the finite family of ideals $I_j\delta_{(j,j)}$ of $R_j\delta_{(j,j)}$ as in Example \ref{ex-symmetric}, that is, $I_1$ is an ideal of $R_1$ and $I_j=\theta_j(I_1)$. We saw in Example \ref{ex-symmetric} that $R$ is $\mcI$-symmetric. Then, by Corollary \ref{cor:ideals}, we have that the subset $I$ of $R$ given in \eqref{expl-ideal} is an ideal. Notice that $I=(I_{jk})_{j,k\in \bI_n}$, where $\theta_j(I_1)\delta_{(k,j)}$.}
\end{exe}
\section{Partial actions of groups on generalized matrix rings}\label{pgma} In what follows in this section, $\tG$ is a group and $R=(M_{ij})_{i,j\in \bI_n}$ is the generalized matrix ring given in \eqref{def-gmr}. We will assume that $R_i=M_{ii}$ is a unital ring, for all $i\in \bI_n$. The main purpose here is to construct a partial action of $\tG$ on $R$ induced by partial actions of $\tG$ on $R_i$, $i\in \bI_n$. 
\subsection{A partial action on  $R$  }\label{sec-4.1}
Assume that for each $i\in \bI_n$, we have a partial action 
\begin{align}\label{par-action-ri}
 \af^{(i)}=\left(D^{(i)}_g,\af^{(i)}_g\right)_{g\in \tG}
\end{align}
of $\tG$ on $R_i$. For each $g\in \tG$, we consider $\mcI_g=\big\{D^{(i)}_g\,:\,i\in\bI_n\big\}$. Suppose that $R$ is $\mcI_g$-symmetric, for all $g\in \tG$. Then
\begin{align}\label{condi-ideals}
	D^{(i)}_gM_{ij}=M_{ij}D^{(j)}_g, \text{ for all } g\in \tG, \,\,i,j\in\bI_n.
\end{align}
By Corollary \ref{cor:ideals}, the subset $I_g$ of $R$, given by 
\begin{align}\label{def-ideal}
I_{g}=\left(D^{(i)}_gM_{ij}\right)_{i,j\in \bI_n},
\end{align}
is an ideal of $R$, for each $g\in \tG$. In order to construct a ring isomorphism between $I_{g\m}$ and $I_g$, we will assume that there exists a finite set of additive bijections 
\begin{align}
	\label{bijec-gamma} \gamma_g=\left\{\gamma^{(ij)}_g:D^{(i)}_{g\m}M_{ij}\to D^{(i)}_gM_{ij}\right\}_{i,j\in \bI_n},\,\,g\in \tG.
\end{align}
Finally, we also assume that: for all $g,h\in \tG$, $i,j,k\in \bI_n$, 
\begin{align}
	\label{cond1} &\gamma^{(ii)}_g=\alpha^{(i)}_g,\quad\quad\gamma^{(ij)}_e=\id_{M_{ij}},& \\[.3em]
	\label{cond2} &\gamma^{(ik)}_g(u)\gamma^{(kj)}_g(v)=\gamma^{(ij)}_g(uv),\,\,\,\text{for all } u\in D^{(i)}_{g\m}M_{ik},\,v\in D^{(k)}_{g\m}M_{kj},\\[.3em]
	\label{cond03} &\big(\gamma^{(ij)}_h\big)^{\m}\!\!\left(D^{(i)}_{g\m}M_{ij}\cap D^{(i)}_{h}M_{ij}\right) \subset D^{(i)}_{(gh)\m}M_{ij},& \\[.3em]
	\label{cond4} &\gamma^{(ij)}_g\left(\gamma^{(ij)}_h(a)\right)=\gamma^{(ij)}_{gh}(a),\,\,\,\text{for all }  a\in \big(\gamma^{(ij)}_h\big)^{\m}\!\!\left(D^{(i)}_{g\m}M_{ij}\cap D^{(i)}_{h}M_{ij}\right).&
\end{align}
Since $\gamma_e^{(ij)}=\id_{M_{ij}}$, it follows from \eqref{cond4} that

\begin{align}
	\label{cond0} &\big(\gamma^{(ij)}_{g}\big)\m= \gamma^{(ij)}_{g\m}, \qquad g\in \tG.
\end{align}

\begin{remark}\label{r1}{\rm Notice that conditions \eqref{cond03} and  \eqref{cond4}   imply that 
		\begin{align}\label{cond5}\big(\gamma^{(ij)}_h\big)^{\m}\!\!\left(D^{(i)}_{g\m}M_{ij}\cap D^{(i)}_{h}M_{ij}\right)=D^{(i)}_{(gh)\m}M_{ij}\cap D^{(i)}_{h\m}M_{ij}, \quad g,h\in \tG,\,\,	 i,j\in \bI_n.
		\end{align}
Indeed, it is clear that $\big(\gamma^{(ij)}_h\big)^{\m}\!\!\left(D^{(i)}_{g\m}M_{ij}\cap D^{(i)}_{h}M_{ij}\right)\subset D^{(i)}_{(gh)\m}M_{ij}\cap D^{(i)}_{h\m}M_{ij},$  for any $g,h\in \tG$. Then, replacing $h$ by $h\m$ and $g$ by $gh$, it follows that
\[\big(\gamma^{(ij)}_{h\m}\big)^{\m}\!\!\left(D^{(i)}_{(gh)\m}M_{ij}\cap D^{(i)}_{h\m}M_{ij}\right)\subset D^{(i)}_{g\m}M_{ij}\cap D^{(i)}_{h}M_{ij}.\] 
Applying $\gamma^{(ij)}_{h\m}$ in the previous inclusion, we get that 
\[D^{(i)}_{(gh)\m}M_{ij}\cap D^{(i)}_{h\m}M_{ij}\subset \gamma^{(ij)}_{h\m}\left(D^{(i)}_{g\m}M_{ij}\cap D^{(i)}_{h}M_{ij}\right).\]
By \eqref{cond0}, we have that $\gamma^{(ij)}_{h\m}=\big(\gamma^{(ij)}_{h}\big)\m$ and hence the result follows.
}
\end{remark}

\begin{def1} \label{def-datum}
	{\rm Let $\alpha^{(i)}$ be a partial action of $\tG$ on $R_i$ as in \eqref{par-action-ri}, $\mcI_g=\big\{D^{(i)}_g\,:\,i\in\bI_n\big\}$ and $\gamma_{g}$ a finite set of additive bijections as in \eqref{bijec-gamma}, for all $i\in \bI_n$ and $g\in \tG$. We say that  \[\mcD=\Big(\{\alpha^{(i)}\}_{i\in \bI_n}, \{\gamma_{g}\}_{g\in\tG}\Big)\]
		is a datum for $R$ if the following statements hold:
		\begin{enumerate}[\rm (a)]
			\item $R$ is $\mcI_g$-symmetric, for all $g\in\tG$, \vspace{.15cm}
			\item the collection $\{\gamma_{g}\}_{g\in\tG}$ satisfies the conditions from \eqref{cond1} to \eqref{cond4}.
		\end{enumerate}
}\end{def1}

Let $\mcD=\big(\{\alpha^{(i)}\}_{i\in \bI_n}, \{\gamma_{g}\}_{g\in\tG}\big)$ be a datum for $R$ as in Definition \ref{def-datum}. For each $g\in \tG$, consider the ideal $I_{g}=\big(D^{(i)}_gM_{ij}\big)_{i,j\in \bI_n}$ of $R$ given by \eqref{def-ideal}. We also define  the additive bijection $\gamma_g:I_{g\m}\to I_g$  by 
\begin{align}\label{bije-pa} 
	\gamma_g=\left(\gamma^{(ij)}_{g}\right)_{i,j\in \bI_n}.
\end{align}
Explicitly, if $u=(a_{ij})_{i,j\in \bI_n}$ is an element of $I_{g\m}$, then  $\gamma_g(u)=\big(\gamma^{(ij)}_{g}(a_{ij})\big)_{i,j\in \bI_n}$.

\begin{teo1}\label{teo-par} 
Let $\mcD=\big(\{\alpha^{(i)}\}_{i\in \bI_n}, \{\gamma_{g}\}_{g\in\tG}\big)$ be a datum for $R$, $I_g$ the ideal of $R$ given in \eqref{def-ideal} and $\gamma_{g}$ a finite set of additive bijections as in \eqref{bije-pa}, for all $g\in \tG$. Then the family of pairs \[\gamma=\left(I_g,\gamma_g\right)_{g\in \tG}\]
is a partial action of the group $\tG$ on the generalized matrix ring $R$. 
\end{teo1}

\begin{proof}
Let $i,j\in \bI_n$ and $g\in \tG$. By \eqref{def-ideal}, we have that \[I^{(ij)}_e=D^{(i)}_eM_{ij}=R_iM_{ij}=M_{ij},\] and consequently $I_e=R$. It follows from 	\eqref{cond1} that $\gamma_e=\id_{R}$. Since $\gamma_g^{(ij)}$ is a bijection, for all $i,j\in\bI_n$, we conclude that $\gamma_g$ is a bijection. Now, we prove that $\gamma_g$ is a ring isomorphism. Let  $a=\big(a^{(ij)}_{g\m}\big)$ and $b=\big(b^{(ij)}_{g\m}\big)$ in $I_{g\m}$, that is, $a^{(ij)}_{g\m},\,b^{(ij)}_{g\m}\in D^{(i)}_{g\m}M_{ij}$, for all $i,j\in \bI_n$. Observe that $\gamma_g(a)\gamma_g(b)=\big(\gamma_g^{(ij)}\big(a^{(ij)}_{g\m}\big)\big)\big(\gamma_g^{(ij)}\big(b^{(ij)}_{g\m}\big)\big)=\big(d^{(ij)}_{g}\big)$, where
\begin{align*}
	d^{(ij)}_{g}=&\,\,\gamma_g^{(i1)}\Big(a^{(i1)}_{g\m}\Big)\gamma_g^{(1j)}\Big(b^{(1j)}_{g\m}\Big)+\ldots+\gamma_g^{(in)}\Big(a^{(in)}_{g\m}\Big)\gamma_g^{(nj)}\Big(b^{(nj)}_{g\m}
	\Big)\\[.3em]
\overset{\mathclap{\eqref{cond2}}}{=}&\,\,\gamma_g^{(ij)}\Big(a^{(i1)}_{g\m}b^{(1j)}_{g\m}\Big)+\ldots+\gamma_g^{(ij)}\Big(a^{(in)}_{g\m}b^{(nj)}_{g\m}\Big)\\[.4em]
                                     =&\,\, \gamma_g^{(ij)}\Big(a^{(i1)}_{g\m}b^{(1j)}_{g\m}+\ldots+ a^{(in)}_{g\m}b^{(nj)}_{g\m}\Big).
\end{align*} 
Since $ab=\big(c^{(ij)}_{g}\big)$ with $c^{(ij)}_{g}=a^{(i1)}_{g\m}b^{(1j)}_{g\m}+\ldots+ a^{(in)}_{g\m}b^{(nj)}_{g\m}$, it follows that $\gamma_g(ab)=\gamma_g(a)\gamma_g(b)$ and hence $\gamma_g$ is a ring isomorphism. Thus, notice that \eqref{cond1}, \eqref{cond03}  and \eqref{cond4} imply that 
$\gamma=\left(I_g,\gamma_g\right)_{g\in \tG}$ is a partial action of $\tG$ on $R$.  
\end{proof}
                    
Let $\mcD=\big(\{\alpha^{(i)}\}_{i\in \bI_n}, \{\gamma_{g}\}_{g\in\tG}\big)$ be a datum for $R$, $I_g$ the ideal of $R$ given in \eqref{def-ideal} and $\gamma_{g}$ a finite set of additive bijections as in \eqref{bije-pa}, for all $g\in \tG$.
For each $k\in \bI_n$, consider the natural inclusions $\iota_k:R_k\to R$ given by $\iota_k(a)=(a_{ij})$, where $a_{kk}=a$ and $a_{ij}=0$ if $i\neq j$. Let $\gamma=\left(I_g,\gamma_g\right)_{g\in \tG}$ be the partial action of $\tG$ on the generalized matrix ring $R$ given in Theorem \ref{teo-par}. Denote by $\gamma_k$ the partial action of $\tG$ on $\iota_k(R_k)$ obtained by restriction, that is, 
\begin{align}\label{part-gamma}
	\gamma^{(k)}=\left(\iota_k\big(D^{(k)}_g\big),\gamma^{(k)}_g\right)_{g\in \tG},
\end{align}
where $\gamma^{(k)}_g=\gamma_g|_{\iota_k(D^{(k)}_{g\m})}$.

\begin{cor1}\label{cor:restric}
	Let $\mcD=\big(\{\alpha^{(i)}\}_{i\in \bI_n}, \{\gamma_{g}\}_{g\in\tG}\big)$ be a datum for $R,$  $\gamma^{(k)}$ be the partial action of $\tG$ on $\iota_k(R_k)$ defined in \eqref{part-gamma} and   $\gamma=\left(I_g,\gamma_g\right)_{g\in \tG}$ be the partial action of $\tG$ on  $R$ given in Theorem \ref{teo-par}. Then $\alpha^{(k)}$ and $\gamma^{(k)}$ are  equivalent partial group actions and $\gamma$ is an extension of $\gamma^{(k)}$, for all $k\in \bI_n$.
\end{cor1}
\begin{proof}
	Let $k\in \bI_n$ and $g\in \tG.$ Notice that $\iota_k:D^{(k)}_g\to \iota_k(D^{(k)}_g)$ is a ring isomorphism and $\iota_k\alpha^{(k)}_g=\gamma^{(k)}_g\iota_k$, which means that  $\alpha^{(k)}$ and $\gamma^{(k)}$ are equivalent. It is immediate that  $\gamma$ is an extension of $\gamma^{(k)}$. 
\end{proof}

\subsection{Conditions for $\gamma$ to be unital}\label{unit} In this subsection,  $\mcD=\Big(\{\alpha^{(i)}\}_{i\in \bI_n}, \{\gamma_{g}\}_{g\in\tG}\Big)$ will denote a datum for $R$ in the sense of Definition \ref{def-datum}. Using Theorem \ref{teo-par}, we will consider the partial action  $\gamma=\left(I_g,\gamma_g\right)_{g\in \tG}$ of $\tG$ on $R$, where the ideal $I_g$ of $R$ and the algebra isomorphism $\gamma_g$ are given respectively by
\[I_{g}=\big(D^{(i)}_gM_{ij}\big)_{i,j\in \bI_n},\qquad\gamma_g=\left(\gamma^{(ij)}_{g}\right)_{i,j\in \bI_n}.\]
We also assume that $\alpha^{(i)}$ is unital, for all $i\in \bI_n$. We will denote by $1^{(i)}_g$ the identity of $D^{(i)}_{g}$, that is, $D^{(i)}_{g}=R_i1^{(i)}_g$.
The purpose of this subsection is to give sufficient condition for $\gamma$ to be a unital partial action. 

Let $g\in \tG$ and consider the diagonal matrix 
\begin{align}\label{unity-for-gamma}
	1_g:=\operatorname{diag}\big(1^{(1)}_g,\ldots,1^{(n)}_g\big)=  \left( \begin{matrix} 1^{(1)}_{g} & 0 & \ldots & 0 \\ 
		0 & 1^{(2)}_{g} & \ddots & \vdots  \\ 
		\vdots & \ddots & \ddots & 0 \\ 
		0 &\ldots & 0 & 1^{(n)}_{g} \end{matrix} \right) \in R.
\end{align}
It is clear that $1_g^2=1_g$, that is, $1_g$ is an idempotent of $R$. Hence, $1_g$ is a natural candidate to be the identity of $I_g$. 
In this direction, the next result gives us sufficient conditions for $1_g$ to be a central element of $R$.

\begin{lem1}\label{lem-central-idempotente}
	Let $g\in \tG$ and let $1_g=\operatorname{diag}\big(1^{(1)}_g,\ldots,1^{(n)}_g\big)$ be the diagonal matrix given in \eqref{unity-for-gamma}. Then $1_g$ is a central element of $R$ if and only if $1^{(i)}_gm_{ij}=m_{ij}1^{(j)}_g$ for all $i,j\in \bI_n$, $m_{ij}\in M_{ij}$ and $g\in \tG$.
\end{lem1}

\begin{proof}
	Let $g\in \tG$, $i,j\in \bI_n$ and $r=(m_{ij})_{i,j\in \bI_n}\in R$. Since $1^{(i)}_gm_{ij}=m_{ij}1^{(j)}_g$, it follows that $1_gr=r1_g$ and hence $1_g$ is a central element of $R$. Conversely, suppose that  $1_g$ is a central element of $R$. Let  $i,j\in \bI_n$ and $m_{ij}\in M_{ij}$. Now, take  $r_{i,j}=(m'_{kl})_{k,l\in \bI_n}\in R$, where $m'_{kl}=0$ for $(k,l)\neq (i,j)$ and $m'_{ij}=m_{ij}$. Since $1_gr_{i,j}=r_{i,j}1_g$, it follows that  $1^{(i)}_gm_{ij}=m_{ij}1^{(j)}_g$. 
\end{proof}

\begin{remark}\label{obs-igualdade}
	{\rm Let $g\in \tG$ and $i,j\in \bI_n$.  Since  $1^{(i)}_gM_{ij}\subset D^{(i)}_gM_{ij}=1^{(i)}_gR_iM_{ij}\subseteq 1^{(i)}_gM_{ij}$, 
	we obtain that $D^{(i)}_gM_{ij}=1^{(i)}_gM_{ij}$.
	}
	
\end{remark}

\begin{lem1}\label{lem-unital-action-gamma} For each $g\in \tG$, let  $1_g=\operatorname{diag}\big(1^{(1)}_g,\ldots,1^{(n)}_g\big)$ be the diagonal matrix given in \eqref{unity-for-gamma}. If  $1^{(i)}_gm_{ij}=m_{ij}1^{(j)}_g$, for all $i,j\in \bI_n$ and $m_{ij}\in M_{ij}$, then $I_g=R1_g$.
\end{lem1}
\begin{proof}
	Let $g\in \tG$, $i,j\in \bI_n$ and $I_{ij}=D^{(i)}_{g}M_{ij}$. By Lemma \ref{lem-central-idempotente}, $1_g$ is a central element of $R$. Now, notice that $I_g=R1_g$ if and only if 
	
	\[ \left( \begin{matrix} D^{(1)}_{g} & I_{12} & \ldots & I_{1n} \\ 
		I_{21} & D^{(2)}_{g} & \ddots & \vdots  \\ 
		\vdots & \ddots & \ddots & I_{(n-1)n} \\ 
		I_{n1} &\ldots & I_{n(n-1)} & D^{(n)}_{g} \end{matrix} \right)=\left( \begin{matrix} D^{(1)}_{g} & M_{12}1^{(2)}_g & \ldots & M_{1n}1^{(n)}_g \\ 
		M_{21}1^{(1)}_g & D^{(2)}_{g} & \ddots & \vdots  \\ 
		\vdots & \ddots & \ddots & M_{(n-1)n}1^{(n)}_g \\ 
		M_{n1}1^{(1)}_g &\ldots & 	M_{n1}1^{(n-1)}_g & D^{(n)}_{g} \end{matrix} \right). \]
	By Remark \ref{obs-igualdade}, $I_{ij}=D^{(i)}_{g}M_{ij}=1^{(i)}_gM_{ij}=M_{ij}1^{(j)}_g$, and consequently $I_g=R1_g$.
\end{proof}

\begin{prop1}\label{prop-gamma-unital}
For each $g\in \tG$, let  $1_g=\operatorname{diag}\big(1^{(1)}_g,\ldots,1^{(n)}_g\big)$ be the diagonal matrix given in \eqref{unity-for-gamma}. If  $1^{(i)}_gm_{ij}=m_{ij}1^{(j)}_g$, for all $i,j\in \bI_n$ and $m_{ij}\in M_{ij}$, then $\gamma$ is a unital partial action of $\tG$ on $R$.	
\end{prop1}
\begin{proof}
	The result follows directly from Lemmas \ref{lem-central-idempotente} and \ref{lem-unital-action-gamma}.
\end{proof}

\subsection{Induced partial actions}\label{subsec-4-3} Let $\beta^{(i)}$ be a global action of $\tG$ on $R_i$, for all $i\in \bI_n$. Let $g\in\tG$ and $\mcI_g=\{R_i\,:\,i\in \bI_n\}$, that is, we are taking $\mcI_g=\mcI_h$ for all $g,h\in \tG$. Since $R^{(k)}$ is a unital ring, $k\in \bI_n$,  we have that $R_iM_{ij}=M_{ij}=M_{ij}R_j$. Thus $R$ is $\mcI_g$-symmetric. We also assume that there exists a finite set of additive bijections 
\begin{align*}
	\beta_g=\left\{\beta^{(ij)}_g:M_{ij}\to M_{ij}\right\}_{i,j\in \bI_n}, \quad i,j\in\bI_n,\,\,g\in \tG,
\end{align*}
such that: for all  $g,h\in \tG$, $i,j,k\in\bI_n$,
\begin{align}
	\label{glo-con1}	&\beta^{(ii)}_g=\beta^{(i)}_g,\quad\quad\beta^{(ij)}_e=\id_{M_{ij}}, \qquad \beta^{(ij)}_g\circ\beta^{(ij)}_h=\beta^{(ij)}_{gh},& \\[.3em]
	\label{glo-con2}	&\beta^{(ik)}_g(u)\beta^{(kj)}_g(v)=\beta^{(ij)}_g(uv),\,\,\,\text{for all } u\in M_{ik},\,v\in M_{kj}.
\end{align}
Under these conditions, it is clear that $\tilde{\mcD}=\Big(\{\beta^{(i)}\}_{i\in \bI_n}, \{\beta_{g}\}_{g\in\tG}\Big)$ is a datum for $R$ which will be called {\it a global datum} for $R$.
Given $r=(m_{ij})_{i,j\in \bI_n}\in R$ and $g\in \tG$, we define the ring isomorphism $\beta_g:R\to R$ given by
\[\beta_{g}(r)=\big(\beta^{(ij)}_{g}(m_{ij})\big)_{i,j\in \bI_n}.\] 
By using Theorem \ref{teo-par}, one would get that $\beta=\{\beta_{g}\}_{g\in \tG}$ is a global action of $\tG$ on $R$.\smallbreak

Now consider $\mathcal{J}=\{J_k\subset R_k\,:\,J_k \text{ is an ideal of }R_k,\text{ for all }k\in \bI_n\}$ a finite family of ideals and assume that $R$ is $\mathcal{J}$-symmetric, that is, $J_kM_{kl}=M_{kl}J_l$, for all $k,l\in \bI_n$. Then, it follows from Propositions \ref{prop-symmetric} and \ref{prop:ideals-gm} that 
\begin{align}\label{ring-J}
	J=\big(J_{ij}\big)_{i,j\in \bI_n},\,\text{ where } J_{ij}=J_iM_{ij},
\end{align}
is an ideal of $R$. Hence, $\beta$ induces a partial action $\alpha=(\alpha_g,J_g)$ of $\tG$ on the generalized matrix ring $J$ given by
\begin{align}\label{pa-induced}
	J_g=J\cap \beta_{g}(J) ,\qquad\quad \alpha_{g}=\beta_g|_{J_{g\m}}:J_{g\m}\to J_{g}.
\end{align}
Notice also that, for each $i\in \bI_n$, we have a partial action $\alpha^{(i)}=\big(\alpha^{(i)}_g,D^{(i)}_g\big)$ of $\tG$ on $J_i$ which is induced by $\beta^{(i)}$, that is, 
\begin{align}\label{pa-Ji}
	D^{(i)}_g=J_i\cap \beta^{(i)}_g(J_i),\qquad\quad \alpha^{(i)}_g=\beta^{(i)}_g|_{D^{(i)}_{g\m}}:D^{(i)}_{g\m}\to D^{(i)}_{g}.	
\end{align}
Since $\alpha^{(i)}$ is a partial action of $\tG$ on $J_i$, $i\in \bI_n$, we have the following natural questions:

\begin{enumerate}[\,\,\rm (a)]
	\item Is it possible to construct the partial action $\gamma$ of $\tG$ on $J$ as in $\S\,$\ref{sec-4.1}?\smallbreak
	\item Suppose that the answer to question (a) is affirmative. In this case, do the partial actions $\gamma$ and $\alpha$ of $\tG$ on $J$ coincide? 
\end{enumerate}

In order to answer the previous questions, in what follows in this subsection, we will assume  that $J_i=R_ie_i$, where $e_i$ is a central idempotent of $R_i.$
Then, 
\begin{align}\label{aux-flat-1}
D^{(i)}_{g}:=J_i\cap \beta^{(i)}_g(J_i)=R_ie_i\beta^{(i)}_g(e_i)=J_i\beta^{(i)}_g(J_i). 
\end{align}
Thus, \[D^{(i)}_{g}J_{ij}=J_i\beta^{(i)}_g(J_i)J_iM_{ij}=J_i\beta^{(i)}_g(J_i)M_{ij}=D^{(i)}_{g}M_{ij}.\]
Moreover, we have that $D^{(i)}_{g\m}\otimes_{R_i }M_{ij}\simeq D^{(i)}_{g\m}M_{ij}$ and 
$D^{(i)}_{g}\otimes_{R_i }M_{ij}\simeq D^{(i)}_{g}M_{ij}$ as additive groups. 
On the other hand, we recall that $\beta^{(ij)}_g:M_{ij}\to M_{ij}$ is an additive bijection and $\alpha^{(i)}_g:D^{(i)}_{g\m}\to D^{(i)}_{g}$ is a ring isomorphism. Then, we have that 
$\gamma^{(ij)}_g:D^{(i)}_{g\m}M_{ij}\to D^{(i)}_{g}M_{ij}$ defined by
\begin{align}\label{aux-flat-2}
\gamma^{(ij)}_g\big(e_i\beta^{(i)}_{g\m}(e_i)m\big)=e_i\beta^{(i)}_g(e_i)\beta^{(ij)}_{g}m,\qquad m\in M_{ij},\,\,i,j\in \bI_n,
\end{align}
is an addictive bijection. Hence, for each $g\in \tG$, we can consider the family of additive bijections
\begin{align}\label{gamma-specific}
	\gamma_g=\left\{\gamma^{(ij)}_g:D^{(i)}_{g\m}M_{ij}\to D^{(i)}_gM_{ij}\right\}_{i,j\in \bI_n}.
\end{align}

\begin{prop1}\label{coincide} 
Let $\alpha^{(i)}$  be the partial action of $\tG$ on $J_i$ given  by \eqref{pa-Ji} and $\gamma_{g}$ the family of additive bijections defined in \eqref{aux-flat-2} and \eqref{gamma-specific}.	Then the pair
\[\mcD=\Big(\{\alpha^{(i)}\}_{i\in \bI_n}, \{\gamma_{g}\}_{g\in\tG}\Big)\]
is a datum for $J$. Moreover, the correspondent partial action $\gamma$ of $\tG$ on $J$, associated to $\mcD$ as in $\S\,$\ref{sec-4.1}, coincides with the partial action $\alpha$ induced by $\beta$ and defined in \eqref{pa-induced}.
\end{prop1}
\begin{proof}
	Let $g\in \tG$ and $i,j\in \bI_n$. Observe that
\begin{align*}
	D^{(i)}_{g}J_{ij}&=D^{(i)}_{g}M_{ij}=J_i\beta^{(i)}_g(J_i)M_{ij}=J_i\beta^{(ii)}_g(J_i)\beta^{(ij)}_{g}(M_{ij})\stackrel{\eqref{glo-con2}}{=}J_i\beta^{(ij)}_g(J_iM_{ij})\\[.1em]
		               &=J_i\beta^{(ij)}_g(M_{ij}J_j)\stackrel{\eqref{glo-con2}}{=}J_i\beta^{(ij)}_g(M_{ij})\beta^{(j)}_g(J_j)=J_iM_{ij}\beta^{(j)}_g(J_j)=M_{ij}J_j\beta^{(j)}_g(J_j)\\[.12em]
		               &=J_{ij}D^{(j)}_{g},
\end{align*} 
which implies that $J$ is $\mathcal{I}$-symmetric, where $\mathcal{I}=\big\{D^{(i)}_g\,:\,i\in\bI_n\big\}$. It is straightforward to check that the collection $\{\gamma_{g}\}_{g\in\tG}$ satisfies the conditions from \eqref{cond1} to \eqref{cond4} and hence $\mcD$ is a datum for $J$. We recall from $\S$ \ref{sec-4.1} that the corresponding partial action of $\tG$ on $J$ associated to $\mcD$ is $\gamma=\left(I_g,\gamma_g\right)_{g\in \tG}$, where  $I_{g}=\big(D^{(i)}_gJ_{ij}\big)_{i,j\in \bI_n}=\big(D^{(i)}_gM_{ij}\big)_{i,j\in \bI_n}$ and $\gamma_g=\big(\gamma^{(ij)}_{g}\big)_{i,j\in \bI_n}$. Then, we need to prove that $D^{(i)}_gM_{ij}=J_iM_{ij}\cap \beta^{(ij)}_{g}(J_iM_{ij})$ and that $\gamma_{g}=\alpha_{g}$, for all $g\in \tG$ and $i,j\in \bI_n$. First, it is immediate that \[D^{(i)}_gM_{ij}=J_i\beta^{(i)}_{g}(J_i)M_{ij}\subset J_iM_{ij}\cap\beta^{(i)}_{g}(J_i)\beta^{(ij)}_{g}(M_{ij})=J_iM_{ij}\cap \beta^{(ij)}_{g}(J_iM_{ij}).\] Now, let $x\in J_iM_{ij}\cap\beta^{(i)}_{g}(J_i)\beta^{(ij)}_{g}(M_{ij})$. Then, there are elements $m,\,\tilde{m}\in M_{ij}$ such that $x=e_im=\beta^{(i)}_{g}(e_i)\tilde{m}$. Thus, 
\[x=e_im=e_i(e_im)=e_ix=e_i\beta^{(i)}_{g}(e_i)\tilde{m}\in e_i\beta^{(i)}_{g}(e_i)M_{ij}=D^{(i)}_gM_{ij},\]
and consequently $D^{(i)}_gM_{ij}=J_iM_{ij}\cap \beta^{(ij)}_{g}(J_iM_{ij})$. Finally, given $y\in D^{(i)}_{g\m}M_{ij}$, there exists $m\in M_{ij}$ such that $y=e_i\beta^{(i)}_{g\m}(e_i)m$ and we have that
\begin{align*}
\gamma^{(ij)}_g(y)&=e_i\beta^{(i)}_g(e_i)\beta^{(ij)}_{g}(m)=\beta^{(i)}_g\big( e_i\beta^{(i)}_{g\m}(e_i)\big)\beta^{(ij)}_{g}(m)\\
                  &\stackrel{\eqref{glo-con2}}{=}\beta^{(ij)}_g\big( e_i\beta^{(i)}_{g\m}(e_i)m\big)=\beta^{(ij)}_g(y).
\end{align*}
Hence $\gamma^{(ij)}_g=\beta^{(ij)}_g$ in $D^{(i)}_{g\m}M_{ij}$, for all $i,j\in \bI_n$ and $g\in \tG$.
Moreover, $\alpha_{g}=\beta_g|_{J_{g\m}}$ and \[J_{g\m}=J\cap \beta_{g\m}(J)=\Big(J_iM_{ij}\cap\beta^{(ij)}_{g\m}(J_iM_{ij})\Big)_{i,j\in \bI_n}=\Big(D^{(i)}_{g\m}M_{ij}\Big)_{i,j\in \bI_n}.\]
Therefore, $\alpha=\gamma$.
\end{proof}

\subsection{Morita equivalent partial actions}\label{sec-exis-gamma}
Throughout this subsection, $\tG$ is a group and $R=(M_{ij})_{i,j\in \bI_n}$ is a generalized matrix ring as in \eqref{def-gmr} such that $R_i=M_{ii}$ is a unital ring, for all $i\in \bI_n$. 
Also, $k$ denotes a  unital  commutative ring and $\alpha=\left(D_g,\af_g\right)_{g\in \tG}$ is a partial action of $\tG$ on a $k$-algebra $A$. 

Our aim in this section is to investigate the relation between the notions of datum and Morita equivalence for partial group actions.

\begin{def1} \label{amod}{\rm (\cite[Def. 2.10]{ADES}) Assume that $\alpha$ is regular. A left unital $A$-module $M$ together with a family of $k$-module isomorphisms $\{\gamma_{g}:D_{g\m}M\to D_gM\}_{g\in \tG}$ is a {\it left $\af$-module} if the following properties are satisfied: for all $g,h\in \tG$,
\begin{enumerate}[\rm (i)]
	\item $\gamma_e$ is the identity map of $M$,\vspace{.1cm}
	\item $\gamma_{g}\gamma_{h}(m)=\gamma_{gh}(m)$, for all $m\in D_{h\m}D_{(gh)\m}M$,\vspace{.1cm}
	\item $\gamma_g(am)=\alpha_g(a)\gamma_g(m)$, for all $a\in D_{g\m}$ and $m\in D_{g\m}M$.
	\end{enumerate}}	
\end{def1}
Similarly one could define the notion of a right $\alpha$-module. Moreover, if $\alpha'$ is a regular partial action of $\tG$ on a $k$-algebra $A'$, then by an $(\alpha,\alpha')$-bimodule  we understand a unital $(A,A')$-bimodule which is a left $\alpha$-module and a right $\alpha'$-module.

\begin{prop1}\label{prop-implica} 
Let $\mcD=\big(\{\alpha^{(i)}\}_{i\in \bI_n}, \{\gamma_{g}\}_{g\in\tG}\big)$ be a datum for $R$ and $i,j\in \bI_n$. Assume that $\alpha^{(i)},\,\alpha^{(j)}$ are regular and $R_i,\,R_j$ are $k$-algebras. Then:
\begin{enumerate}[\rm (i)]
\item the left $R_i$-module $M_{ij}$ together with the family $\{\gamma^{(ij)}_g:D^{(i)}_{g\m}M_{ij}\to D^{(i)}_gM_{ij}\}_{g\in \tG}$  of $k$-module isomorphisms is a left $\alpha^{(i)}$-module. 
\item the right $R_j$-module $M_{ij}$ together with the family $\{\gamma^{(ij)}_g:M_{ij}D^{(j)}_{g\m}\to M_{ij} D^{(j)}_g\}_{g\in \tG}$  of $k$-module isomorphisms is a right $\alpha^{(j)}$-module. 
\end{enumerate}
\end{prop1}
\begin{proof}
Notice that \eqref{cond1} implies that $\gamma^{(ij)}_e$ is the identity map of $M_{ij}$. Also, if  $k=i$ in \eqref{cond2} then we get by  \eqref{cond1} and  \eqref{cond2}  that 
$\gamma^{(ij)}_g(am)=\alpha^{(i)}_g(a)\gamma^{(ij)}_g(m)$, for all $a\in D^{(i)}_{g\m}$ and $m\in D^{(i)}_{g\m}M_{ij}$. Finally, since $\alpha^{(i)}$ is regular, we have that 
\[D^{(i)}_{h\m}D^{(i)}_{(gh)\m}M_{ij}=\big(D^{(i)}_{h\m}\cap D^{(i)}_{(gh)\m}\big)M_{ij}\subset D^{(i)}_{h\m} M_{ij}\cap D^{(i)}_{(gh)\m}M_{ij}.\]
Then, given $m\in D^{(i)}_{h\m}D^{(i)}_{(gh)\m}M_{ij}$, it follows that
\[m\in D^{(i)}_{h\m} M_{ij}\cap D^{(i)}_{(gh)\m}M_{ij}\stackrel{\eqref{cond5}}{=}\big(\gamma^{(ij)}_h\big)^{\m}\!\!\left(D^{(i)}_{g\m}M_{ij}\cap D^{(i)}_{h}M_{ij}\right).\]
Thus, \eqref{cond4} implies that  $\gamma^{(ij)}_g\gamma^{(ij)}_h(m)=\gamma^{(ij)}_{gh}(m)$ and hence (i) follows. By \eqref{condi-ideals}, it is clear that the proof of (ii) is similar.
\end{proof}

As an immediate consequence of the previous result, we have the following.

\begin{cor1}\label{cor-bimodule}
Let $\mcD=\big(\{\alpha^{(i)}\}_{i\in \bI_n}, \{\gamma_{g}\}_{g\in\tG}\big)$ be a datum for $R$ and $i,j\in \bI_n$. Assume that $\alpha^{(i)},\,\alpha^{(j)}$ are regular and $R_i,\,R_j$ are $k$-algebras. Then $M_{ij}$ is an $(\alpha^{(i)},\alpha^{(j)})$-bimodule.
\end{cor1}

Now, we recall from  \cite[Def. 2.8]{ADES} the notion of Morita equivalence for partial group actions.
\begin{def1}\label{def-Morita-equivalence}{\rm
Let $\alpha=\left(D_g,\af_g\right)_{g\in \tG}$ and $\alpha'=\left(D'_g,\af'_g\right)_{g\in \tG}$ be regular partial actions of $\tG$ on the $k$-algebras $A$ and $A'$, respectively. We say that $\alpha$ and $\alpha'$ are Morita equivalent if the following conditions hold:
\begin{enumerate}[\rm (i)]
	\item there exists a Morita context $(A,A',M,M',\tau,\tau')$ such that $\tau,\tau'$ are surjective, $M,M'$ are unital bimodules and $M'D_gM=D'_g$, for all $g\in \tG$,\vspace{.1cm}
	\item there exists a product partial action $\theta=\left(\bar{D}_g,\bar{\af}_g\right)_{g\in \tG}$ of $\tG$ on the Morita ring
	$$ \left( \begin{matrix} A & M \\ M' & A' \end{matrix} \right)$$
	such that the restriction of $\theta$ to $\left( \begin{matrix} A & 0 \\ 0 & 0 \end{matrix} \right)$ is $\alpha$, whereas the restriction of $\theta$ to $\left( \begin{matrix} 0 & 0 \\ 0 & A' \end{matrix} \right)$ is $\alpha'$. 
\end{enumerate}}
\end{def1}

For the  convenience of the  reader, we recall \cite[Prop. 2.11]{ADES}.
\begin{prop1} \label{cmo}In Definition \ref{def-Morita-equivalence}, item {\rm(ii)} can be replaced by the following:
	there exist an $(\alpha, \alpha')$-bimodule structure on $M$ and an $(\alpha', \alpha)$-bimodule structure on $M'$ such that:
	\begin{enumerate}[\rm (i)]
		\item $\alpha_g(mm')=\gamma_g(m)\gamma'_g(m')$,
		\item $\alpha'_g(m'm)=\gamma'_g(m')\gamma_g(m)$, \vspace{.1cm}
	\end{enumerate}
	for all $m\in D_{g\m}M,\, m'\in D'_{g\m}M'$.
\end{prop1}

In the next two Propositions we investigate the relation between a datum and Morita equivalent partial actions.

\begin{prop1} \label{prop-Morita-1}
Let $\mcD=\big(\{\alpha^{(i)}\}_{i\in \bI_n}, \{\gamma_{g}\}_{g\in\tG}\big)$ be a datum for $R$ and $i,j\in \bI_n$. Assume that $R_i$ and $R_j$ are $k$-algebras, with $M_{ij}M_{ji}=R_i$,  $M_{ji}M_{ij}=R_j$ and that the partial actions $\alpha^{(i)}$ and $\alpha^{(j)}$ are regular. Then $\alpha^{(i)}$ and $\alpha^{(j)}$ are Morita equivalent.
\end{prop1}
\begin{proof}
Notice that $\Upsilon_{ij}=(R_i, R_j, M_{ij}, M_{ji}, \theta_{iji}, \theta_{jij})$ is a Morita context and $\theta_{iji},\, \theta_{jij}$ are surjective. Since $M_{ji}$ is $\mcI_g$-symmetric, it follows that \[M_{ji}D^{(i)}_{g}M_{ij}=D^{(j)}_{g}M_{ji}M_{ij}=D^{(i)}_{g},\,\,\text{ for all }\,g\in\tG.\]
 Hence, $\Upsilon_{ij}$ satisfies (i) of Definition \ref{def-Morita-equivalence}. By Corollary \ref{cor-bimodule}, $M_{ij}$ is an $(\af^{(i)},\af^{(j)})$-bimodule and $M_{ji}$  is an $(\af^{(j)},\af^{(i)})$-bimodule. Let $u\in D^{(i)}_{g\m}M_{ij}$ and $v\in D^{(j)}_{g\m}M_{ji}$. By \eqref{cond1} and \eqref{cond2}, we have that $\alpha^{(i)}_g(uv)=\gamma^{(ii)}_g(uv)=\gamma^{(ij)}_g(u)\gamma^{(ji)}_g(v)$. Similarly, $\alpha^{(j)}_g(vu)=\gamma^{(ji)}_g(v)\gamma^{(ij)}_g(u)$. Thus, it follows from Proposition \ref{cmo} that $\alpha^{(i)}$ and $\alpha^{(j)}$ are Morita equivalent.\end{proof}

\begin{prop1}\label{prop-Morita-2}
Consider regular partial actions $\alpha^{(1)}$ and $\alpha^{(2)}$ of $\tG$ on the $k$-algebras $R_1$ and $R_2$, respectively. Assume that $\alpha^{(1)}$ and $\alpha^{(2)}$ are Morita equivalent and that $(R_1,R_2,M_{12},M_{21},\tau_{12},\tau_{21})$
is a Morita context that satisfies (i) of Definition \ref{def-Morita-equivalence} and 
\begin{align}\label{eq-aux-prop}
	D^{(i)}_gM_{ij}\cap D^{(i)}_hM_{ij}=\big(D^{(i)}_{g}\cap D^{(i)}_{h}\big)M_{ij},\qquad i,j\in \I_2,\,\,g\in\tG.
\end{align}
Then there is a datum
$\mcD=\big(\{\alpha^{(i)}\}_{i\in \bI_2}, \{\gamma_{g}\}_{g\in\tG}\big)$ for the Morita ring
$ R=\left( \begin{matrix} R_1 & M_{12} \\ M_{21} & R_2 \end{matrix} \right).$
\end{prop1}
\begin{proof}
By Proposition \ref{cmo}, there exist an $(\alpha^{(1)}, \alpha^{(2)})$-bimodule structure on $M_{12}$ and $(\alpha^{(2)}, \alpha^{(1)})$-bimodule structure on $M_{21}$. Thus, there are two families of $k$-module isomorphisms
\[\Big\{\gamma^{(12)}_{g}:D^{(1)}_{g\m}M_{12}\to D^{(1)}_gM_{12}\Big\}_{g\in \tG},\quad \left\{\gamma^{(21)}_{g}:D^{(2)}_{g\m}M_{21}\to D^{(2)}_gM_{21}\right\}_{g\in \tG}.\] 
For each $g\in \tG$, we consider the set $\gamma_g=\left\{\gamma^{(11)}_{g}=\alpha^{(1)}_{g},\gamma^{(12)}_{g},\gamma^{(21)}_{g},\gamma^{(22)}_{g}=\alpha^{(2)}_{g}\right\}$ and let $\mcD=\big(\{\alpha^{(i)}\}_{i\in \bI_2}, \{\gamma_{g}\}_{g\in\tG}\big)$. By Definition \ref{def-Morita-equivalence} (i), we have that $M_{21}D^{(1)}_gM_{12}=D^{(2)}_g$, for all $g\in \tG$. It follows from \cite[Rem. 2.9]{ADES}  that
\[D^{(1)}_gM_{12}=M_{12}D^{(2)}_g,\qquad D^{(2)}_gM_{21}=M_{21}D^{(1)}_g,\,\text{ for all } g\in \tG.\]
Thus, $R$ is $\mcI_g$-symmetric, where $\mcI_g=\big\{D^{(1)}_g,D^{(2)}_g\big\}$. Since $M_{12}$ is a lef $\alpha^{(1)}$-module, it follows from Definition \ref{amod} that $\gamma^{(12)}_{e}=\id_{M_{12}}$. Similarly, $\gamma^{(21)}_{e}=\id_{M_{21}}$ and we conclude that \eqref{cond1} is true. Observe that \eqref{cond2} follows from Definition \ref{amod} (i) and Proposition \ref{cmo} (i) and (ii). In order to prove \eqref{cond03}, notice that
\begin{align*}
	\big(\gamma^{(ij)}_h\big)^{\m}\!\!\left(D^{(i)}_{g\m}M_{ij}\cap D^{(i)}_{h}M_{ij}\right)&\overset{\mathclap{\eqref{eq-aux-prop}}}{=} \big(\gamma^{(ij)}_h\big)^{\m}\!\!\left(\big(D^{(i)}_{g\m}\cap D^{(i)}_{h}\big)M_{ij}\right)\\
	&=\big(\gamma^{(ij)}_h\big)^{\m}\!\!\left(D^{(i)}_{g\m}D^{(i)}_{h}M_{ij}\right)\\
	&=\gamma^{(ij)}_{h^{-1}}\!\!\left(D^{(i)}_{g\m}D^{(i)}_{h}M_{ij}\right)\\
	&\overset{(\star)}{\subset} D_{h^{-1}}M_{ij}\cap D_{h^{-1}g^{-1}}M_{ij}\\
	&\overset{\mathclap{\eqref{eq-aux-prop}}}{=} \big( D_{h^{-1}}\cap D_{h^{-1}g^{-1}} \big)M_{ij},
\end{align*}
where the inclusion in $(\star)$ is valid by Definition \ref{amod} (ii). Finally, \eqref{cond4} is immediate from Definition \ref{amod} (ii). Thus, $\mcD$ is a datum for $R$.
\end{proof}

\section{Galois theory}\label{gsep}

Throughout this section, $\tG$ is a finite group and $R=(M_{ij})_{i,j\in \bI_n}$ is the generalized matrix ring given in \eqref{def-gmr}. We will assume that $R_i=M_{ii}$ is a unital ring, for all $i\in \bI_n$. Also,  $\mcD=\Big(\{\alpha^{(i)}\}_{i\in \bI_n}, \{\gamma_{g}\}_{g\in\tG}\Big)$ will denote a datum for $R$ in the sense of Definition \ref{def-datum}. Using Theorem \ref{teo-par}, we will consider the partial action  $\gamma=\left(I_g,\gamma_g\right)_{g\in \tG}$ of $\tG$ on $R$, where the ideal $I_g$ of $R$ and the algebra isomorphism $\gamma_g$ are given respectively by
\[I_{g}=\big(D^{(i)}_gM_{ij}\big)_{i,j\in \bI_n},\qquad\gamma_g=\left(\gamma^{(ij)}_{g}\right)_{i,j\in \bI_n}.\]
We also assume that $\alpha^{(i)}$ is unital, for all $i\in \bI_n$. We will denote by $1^{(i)}_g$ the identity of $D^{(i)}_{g}$, that is, $D^{(i)}_{g}=R_i1^{(i)}_g$. Finally, we will suppose that the central idempotents $1^{(i)}_g$, $i\in \bI_n$ and $g\in \tG$, satisfy the following condition:
\begin{align}\label{cond-identity}
	1^{(i)}_gm_{ij}=m_{ij}1^{(j)}_g, \,\,\text{for all }\,i,j\in \bI_n,\,\,m_{ij}\in M_{ij}.
\end{align}
Then, by Proposition \ref{prop-gamma-unital}, $\gamma$ is a unital action and $I_g=R1_g$, where $1_g$ is the diagonal matrix given in \eqref{unity-for-gamma}.\smallbreak

The purpose in this section is to investigate some properties of the ring extension $R^{\gamma}\subset R$ related to Galois theory. 

\subsection{Invariant subring, trace map and separability}\label{galois}

Let $\alpha=\left(A_g,\alpha_g\right)_{g\in \tG}$ be a unital partial action of $\tG$ on a ring $A$ with $A_g=A\tilde{1}_g$, where $\tilde{1}_g$ is a central idempotent of $A$, for all $g\in \tG$. We recall from \cite{DFP} that the {\it subring of invariants of $A$} is given by 
\[A^{\alpha}=\{a\in A\,:\,\alpha_g(a\tilde{1}_{g\m})=a\tilde{1}_g,\,\,\forall g\in \tG\}.\]
In order to determine the subring of invariants $R^{\gamma}$, we define
\begin{align}\label{invariant-submodule}
	M^{\gamma}_{ij}:=\{m\in M_{ij}\,:\,\gamma^{(ij)}_{g}(m1^{(j)}_{g\m})=m1^{(j)}_g,\,\,\forall g\in \tG\}\subset M_{ij}.
\end{align}
It follows from \eqref{cond1} that $M^{\gamma}_{ii}=R^{\alpha^{(i)}}_i$. Also, we have the following.
\begin{prop1}\label{lem-inv-submodule}
$M^{\gamma}_{ij}$ is an $\big(R^{\alpha^{(i)}}_i,R^{\alpha^{(j)}}_j\big)$-bimodule, for all $i,j\in \bI_n$. 	
\end{prop1}
\begin{proof}
Let $i,j\in \bI_n$, $a\in R^{\alpha^{(i)}}_i$ and $m\in M^{\gamma}_{ij}$. Then, for all $g\in \tG$,
	\begin{align*}
		\gamma^{(ij)}_{g}\big((am)1^{(j)}_{g\m}\big)&=\gamma^{(ij)}_{g}\big(a(m1^{(j)}_{g\m})\big)\overset{\mathclap{\eqref{cond-identity}}}{=}\gamma^{(ij)}_{g}\big(a(1^{(i)}_{g\m}m1^{(j)}_{g\m})\big)\\[.3em]
		&\overset{\mathclap{\eqref{cond2}}}{=} \gamma^{(ii)}_{g}\big(a1^{(i)}_{g\m}\big)\gamma^{(ij)}_{g}\big(m1^{(j)}_{g\m}\big)\overset{\mathclap{\eqref{cond1}}}{=}
		\alpha^{(i)}_{g}\big(a1^{(i)}_{g\m}\big)\gamma^{(ij)}_{g}\big(m1^{(j)}_{g\m}\big)\\[.3em]
		&=a1^{(i)}_{g}m1^{(j)}_g=(am)1^{(j)}_g.
	\end{align*}
Hence $ax\in M^{\gamma}_{ij}$. Similarly, if $b\in R^{\alpha^{(j)}}_j$, then $mb\in M^{\gamma}_{ij}$. Since $M_{ij}$ is an $(R_i,R_j)$-bimodule, it follows that $a(mb)=(am)b$ and the result is proved.
\end{proof}  

\begin{prop1}\label{prop-inva-gamma}
	The subring of invariants $R^{\gamma}$ is given by\vspace{.1cm}
	\[R^{\gamma}=\left( \begin{matrix}  R^{\alpha^{(1)}}_1& M^{\gamma}_{12} & \ldots & M^{\gamma}_{1n} \\[.2em]
		M^{\gamma}_{21} & R^{\alpha^{(2)}}_2 & \ddots & \vdots  \\[.2em]
		\vdots & \ddots & \ddots & M^{\gamma}_{(n-1)n} \\[.2em] 
		M^{\gamma}_{n1} &\ldots & M^{\gamma}_{n(n-1)} & R^{\alpha^{(n)}}_n\end{matrix} \right).\]
\end{prop1}
\begin{proof}
By definition, $r=(m_{ij})_{i,j\in \bI_n}\in R^{\gamma}$ if and only if $\gamma^{(ij)}_{g}\big(m_{ij}1^{(j)}_{g\m}\big)=m_{ij}1^{(j)}_{g}$, for all $g\in \tG$. Thus, $r\in R^{\gamma}$ if and only if $m_{ij}\in M^{\gamma}_{ij}$. Since $M^{\gamma}_{ii}=R^{\alpha^{(i)}}_i$, we obtain the desired result.
\end{proof}

Now, we also recall from \cite{DFP} the notion of trace map. Given a  unital partial action $\alpha=\left(A_g,\alpha_g\right)_{g\in \tG}$ of  $\tG$ on a ring $A$, the {\it trace map associated to $\alpha$} is defined as the map $t_{\alpha}:A\to A^{\alpha}$ given by
\[t_{\alpha}(a)=\sum_{g\in \tG}\alpha_{g}(a\tilde{1}_{g\m}), \,\,\text{ where }\,a\in A\,\text{ and }\,A_g=A\tilde{1}_g.\]
It is proved in Lemma 2.1 of \cite{BFP} that $t_{\alpha}(A)\subset A^{\alpha}$. In order to study the trace map $t_{\gamma}:R\to R^{\gamma}$, we consider the map $t^{(ij)}_{\gamma}:M_{ij}\to M_{ij}$ defined by:
\begin{align}\label{trace-map-aux}
	t^{(ij)}_{\gamma}(m)=\sum_{g\in \tG}\gamma^{(ij)}_g\big(m1^{(j)}_{g\m}\big), \quad \text{for all }\,i,j\in \bI_n,\,\, m\in M_{ij}.
\end{align}
It is straightforward to check that $t^{(ij)}_{\gamma}(M_{ij})\subset M^{\gamma}_{ij}$. Moreover, we have from \eqref{cond1} that $t^{(ii)}_{\gamma}=t_{\alpha^{(i)}}$, for all $i\in \bI_n$.

\begin{lem1}\label{trace-map}
Let $i,j\in \bI_n$ and $t^{(ij)}_{\gamma}$ defined in \eqref{trace-map-aux}. Then the trace map $t_{\gamma}:R\to R^{\gamma}$ is determined by  $t_{\gamma}=\big(t^{(ij)}_{\gamma}\big)_{i,j\in \bI_n}$.
\end{lem1}
\begin{proof}
	Let $r=(m_{ij})_{i,j\in \bI_n}\in R$. Since
	\begin{align*}
		t_{\gamma}(r)&=\sum_{g\in \tG}\gamma_g(r1_{g\m})=\sum_{g\in \tG}\Big(\gamma^{(ij)}_{g}\big(m_{ij}1^{(j)}_{g\m}\big)\Big)_{i,j\in \bI_n}\\[.2em]
		&=\Big(\sum_{g\in \tG}\gamma^{(ij)}_{g}\big(m_{ij}1^{(j)}_{g\m}\big)\Big)_{i,j\in \bI_n}=\Big(t^{(ij)}_{\gamma}(m_{ij})\Big)_{i,j\in \bI_n},
	\end{align*} 
the result is proved.
\end{proof}

Recall that a  unital ring extension $A\subset B$ is called {\it separable} if the multiplication map $m:B\otimes_A B\to B$ is a splitting epimorphism of $A$-bimodules. This is equivalent to saying that there exists an element $e\in B\otimes_AB$ such that $be=eb$, for all $b\in B$, and $m(e)=1_B$. Such an element $e$ is usually called {\it an idempotent of separability} of $B$ over $A$.\smallbreak

Let $\alpha=\left(A_g,\alpha_g\right)_{g\in \tG}$ be a unital partial action of $\tG$ on a ring $A$. It was proved in Theorem 3.1 of \cite{BLP} that $A^{\alpha}\subset A$ is separable if and only if there exists a central element $a\in Z(A)$ such that $t_{\alpha}(a)=1_A$. In order to use this result to characterize when the extension $R^{\gamma}\subset R$ is separable, notice that by a direct calculation we get that
\begin{align} \label{center}
	Z(R)=\left( \begin{matrix} Z(R_1)& 0 & \ldots & 0 \\[.2em]
	0& Z(R_2)& \ddots & \vdots  \\[.2em]
		\vdots & \ddots & \ddots & 0 \\[.2em] 
		0 &\ldots & 0 & Z(R_n)\end{matrix} \right).
\end{align}

\begin{prop1}\label{prop-separability}
	The ring extension $R^{\gamma}\subset  R$ is separable if and only if  the ring extension $R^{\alpha^{(i)}}_i\subset R_i$ is separable, for all $i\in \bI_n$.
\end{prop1}
\begin{proof} Assume that $R^{\gamma}\subset  R$ is separable. Then, it follows from Theorem 3.1 of \cite{BLP} that there is a central element $r=(m_{ij})_{i,j\in \bI_n}\in Z(R)$ such that $t_{\gamma}(r)=1_R$. By \eqref{center}, $m_{ii}\in Z(R_i)$, for all $i\in \bI_n$. Also, Lemma \ref{trace-map} implies that $t_{\alpha^{(i)}}(m_{ii})=1_{R_i}$ and we conclude by Theorem 3.1 of \cite{BLP} that $R^{\alpha^{(i)}}_i\subset R_i$ is separable.
	Conversely, assume that $R^{\alpha^{(i)}}_i\subset R_i$ is separable, for all $i\in \bI_n$. Using again Theorem 3.1 of \cite{BLP}, we obtain $r_i\in Z(R_i)$ such that $t_{\alpha^{(i)}}(r_i)=1_{R_i}$. Thus, by Lemma \ref{trace-map}, the diagonal matrix $r=\operatorname{diag}\big(r_1,\ldots,r_n\big)$ is a central element of $R$ that satisfies $t_{\gamma}(r)=1_R$. Hence, Theorem 3.1 of \cite{BLP} implies that $R^{\gamma}\subset  R$.
\end{proof}

\subsection{Galois extension}
Let $\alpha=\left(A_g,\alpha_g\right)_{g\in \tG}$ be a unital partial action of $\tG$ on a ring $A$ with $A_g=A\tilde{1}_g$, where $\tilde{1}_g$ is a central idempotent of $A$, for all $g\in \tG$. According to \cite{DFP}, the ring extension $A^{\alpha}\subset A$ is {\it a partial Galois extension} whenever there exists a positive integer $m$ and $a_i,b_i\in A$, $i\in \bI_m$, such that \[\sum\limits_{i=1}^ma_i\alpha_g(b_i\tilde{1}_{g\m})=\delta_{e,g}.\]
In this case, the set $\{a_i,b_i\}_{i\in \bI_m}$ is called {\it a partial Galois coordinate system} of $A$ over $A^{\alpha}$. The next result characterizes when $R^{\gamma}\subset R$ is a partial Galois extension.

\begin{teo1}\label{teo-galois}
$R^{\gamma}\subset R$ is a partial Galois extension if and only if $R^{\alpha^{(i)}}_i\subset R_i$ is a partial Galois extension, for all $i\in \bI_n$.	
\end{teo1}
\begin{proof}
Let $i\in \bI_n$ and assume that $\big\{x^{(i)}_j,\,y^{(i)}_j\big\}_{j\in \bI_{m_i}}$ is a partial Galois coordinate system of $R_i$ over $R^{\alpha^{(i)}}_i$, that is,  $\sum\limits_{j=1}^{m_i}x^{(i)}_j\alpha^{(i)}_g\big(y^{(i)}_j1^{(i)}_{g\m}\big)=\delta_{e,g}$. Let $m_{k}:=\max\{m_i\,:\,i\in \bI_n\}$ and for $i\neq k$ we put 
\[x^{(i)}_{m_i+1}=\ldots=x^{(i)}_{m_k}=0,\qquad y^{(i)}_{m_i+1}=\ldots=y^{(i)}_{m_k}=0.\]
Then, for all $j\in \bI_k$, we define the following diagonal matrices
\[r_j:=\operatorname{diag}\big(x^{(1)}_j,\ldots,x^{(n)_j}\big),\qquad s_j:=\operatorname{diag}\big(y^{(1)}_j,\ldots,y^{(n)}_j\big).\]	
Notice  that 
\begin{align*}
	\sum\limits_{j=1}^{m_k}r_j\gamma_g(s_j)&=\operatorname{diag}\Bigg(\sum\limits_{j=1}^{m_1}x^{(1)}_j\alpha^{(1)}_g\big(y^{(1)}_j1^{(1)}_{g\m}\big),\ldots,
	\sum\limits_{j=1}^{m_n}x^{(n)}_j\alpha^{(n)}_g\big(y^{(n)}_j1^{(1)}_{g\m}\big)\Bigg)\\
	&=\operatorname{diag}\big(\delta_{e,g},\ldots,\delta_{e,g}\big)\\
	&=\delta_{e,g}. 
\end{align*}
Thus, $R^{\gamma}\subset R$ is a partial Galois extension.\smallbreak

Conversely, assume that $\big\{r_k=(a_{ij}^k)_{i,j\in \bI_n},\,s_k=(b_{ij}^k)_{i,j\in \bI_n}\big\}_{k\in \bI_l}$ is a partial Galois coordinate system of $R$ over 
$R^{\gamma}$, that is,
	\[\sum\limits_{k=1}^lr_k\gamma_g\big(s_k1_{g\m}\big)=\delta_{e,g},\]
where $1_{g\m}$ is given in \eqref{unity-for-gamma}. Thus, for every $1\leq t\leq n$, we have that
\begin{align}\label{for-galois}
	\sum\limits_{k=1}^l\sum\limits_{j=1}^n a_{tj}^k\gamma^{(jt)}_{g}\big(b_{jt}^k1^{(t)}_{g\m}\big)=\delta_{e,g}=\begin{cases}
		0_{R_t},& \text{ if } g\neq e,\\
		1_{R_t},& \text{ if } g=e.
	\end{cases}	
\end{align}
For all $j\neq t$ and $k\in \bI_l$, let $z_{jk}:=a_{tj}^k\gamma^{(jt)}_{g}\big(b_{jt}^k1^{(t)}_{g\m}\big)\in M_{tj}D^{(j)}_{g}M_{jt}=M_{tj}M_{jt}D^{(t)}_{g}=D^{(t)}_{g}$. Since $\alpha^{(t)}_g$ is a ring isomorphism, there is a unique $y_{jk}\in D^{(t)}_{g\m}$ such that $\alpha^{(t)}_g\big(y_{jk}\big)=z_{jk}$. Now, we consider
\[u_{jk}=\begin{cases}
	1_{R_t},& \text{ if } j\neq t,\\
	a^k_{tt},& \text{ if } j=t,
\end{cases}\qquad \quad v_{jk}=\begin{cases}
	y_{jk},& \text{ if } j\neq t,\\
	b^k_{tt},& \text{ if } j=t,
\end{cases}.\]
Then
\begin{align*}
	\sum\limits_{k=1}^l	\sum\limits_{j=1}^n	u_{jk}\alpha^{(t)}_g\big(v_{jk}1^{(t)}_{g\m}\big)&=\Big(\sum\limits_{k=1}^l\,\,	\sum\limits_{j\in \bI_n,j\neq t} z_{jk}\Big)+ \sum\limits_{k=1}^l	a^k_{tt}\alpha^{(t)}_g\big(b^k_{tt}1^{(t)}_{g\m}\big)\\
	&\overset{\mathclap{\eqref{cond1}}}{=}\Big(\sum\limits_{k=1}^l\,\,	\sum\limits_{j\in \bI_n,j\neq t} z_{jk}\Big)+ \sum\limits_{k=1}^l	a^k_{tt}\gamma^{(tt)}_g\big(b^k_{tt}1^{(t)}_{g\m}\big)\\
	&=\sum\limits_{k=1}^l\sum\limits_{j=1}^n a_{tj}^k\gamma^{(jt)}_{g}\big(b_{jt}^k1^{(t)}_{g\m}\big)\,\,\overset{\mathclap{\eqref{for-galois}}}{=}\,\, \delta_{e,g}.
\end{align*}
Hence, $R^{\alpha^{(t)}}_t\subset R_t$ is a partial Galois extension, for all $t\in \bI_n$. 
\end{proof}

\section{Examples}
In this section we give explicit examples to illustrate the construction of the partial action $\gamma$ of a group $\tG$ on a generalized matrix ring $R$.

\subsection{Partial skew groupoid ring} Throughout this subsection, $\af=(\A_g,\af_g)_{g\in \G}$ denotes a group-type partial action of a connected groupoid $\G$ on the ring $A=\oplus_{y\in \G_0}A_{y}$. We will assume that $\G_0$ is finite and let  $A\star_\af \G$ be the corresponding partial skew groupoid ring. By \eqref{action-beta}, there is a partial action $\beta$ of $\G_0^2$ on $A$. We will see that $A\star_\beta \G_0^2$ is a generalized matrix ring. Also, we will prove that the partial action $\varepsilon$ of the isotropy group $\tG=\G(x)$ on $A\star_\beta \G_0^2$ given in \eqref{pa-oa} is a particular case of the construction presented in the previous subsection. \smallbreak

Notice that $\alpha$ is group-type and hence there exist  $x\in \G_0$ and  a set $\{h_y\,:\,h_y\in \G(x,y)\}$ of morphisms of $\G$ such that $h_x=x$ and $\alpha_{h_y}:A_x\to A_y$, for all $y\in \G_0$.
According to Lemma 4.1 of \cite{BPP}, we have the following global action of $\G_0^2=\G_0\times\G_0$ on $A$:
\begin{align}\label{ac-beta}
	\bt=(A_{(x_i,x_j)},\bt_{(x_i,x_j)})_{(x_i,x_j)\in \G_0^2},\quad A_{(x_i,x_j)}=A_{x_j},\quad \beta_{(x_i,x_j)}=\af_{h_{x_j}}\af_{h\m_{x_i}}.	
\end{align}
Let $A_i=A_{x_i}$, $\theta_{i}=\alpha_{h_{x_i}}$, $A_{(i,j)}=A_j$ and $\theta_{(i,j)}=\theta_{ij}=\theta_j\theta\m_i=\beta_{(x_i,x_j)}$, for all $i,j\in \bI_n$. Then  $\theta=(A_{(i,j)},\theta_{(i,j)})_{(i,j)\in \bI_n^2}$ is  a partial action of $\bI_n^2$ on $B=\prod_{i\in \mathbb{I}_n} A_i$ as in $\S\,$\ref{subsec-skew}. 
\begin{lem1}\label{lem-iso-skews}
	Let $A,B$, $\G$, $\alpha$ and $\theta$ as above. The association $\psi:A*_{\beta}\G_0^2\to B*_{\theta}\bI_n^2$ given by \[\psi(a\delta_{(x_i,y_i)})=a\delta_{(i,j)},\quad a\in A_j,\] is a ring isomorphism. In particular,  $A*_{\beta}\G_0^2$ is a generalized matrix ring.
\end{lem1}
\begin{proof}
	It is straightforward to check that $\psi$ is a ring isomorphism. The last assertion follows from Proposition \ref{prop:skew-gmr}.
\end{proof}

In what follows in this section, $R_i=A_i\delta_{(i,i)}$,  $M_{ij}:=A_{i}\delta_{(j,i)}$, for all $i,j\in \bI_n$, and $R=(M_{ij})_{i,j\in \bI_n}$. Notice that $R$ is a generalized matrix ring with structure given by  \eqref{eq:structure-module-skew} and \eqref{one-more}. \smallbreak

Let $\alpha^{(1)}$ be the partial action of the isotropy group $\tG=\G(x,x)$ on $A_1\delta_{(1,1)}$ obtained by restriction of $\alpha$, that is, $\alpha^{(1)}=\big(\A_g\delta_{(1,1)},\af^{(1)}_g\big)_{g\in \tG}$, where  \[\alpha^{(1)}_g(a_{g\m}\delta_{(1,1)})=\alpha_g(a_{g\m})\delta_{(1,1)}, \quad a_{g\m}\in A_{g\m},\,\,g\in\tG.\]
Then, for each $2\leq i\leq n$, we consider the following partial action of $\tG$ on $A_{i}\delta_{(i,i)}$:
\begin{align} \label{alfa-i-skew}
	\alpha^{(i)}=\big(\af_{h_{x_i}}(A_g)\delta_{(i,i)},\af^{(i)}_g\big)_{g\in \tG},\quad \alpha_g^{(i)}\left(\af_{h_{x_i}}(a_{g\m}) \delta_{(i,i)}\right)=\af_{h_{x_i}}(\alpha_g(a_{g\m}))\delta_{(i,i)},	
\end{align}
for all $a_{g\m}\in A_{g\m}$, $g\in \tG$. We will denote $D^{(i)}_{g}:=\af_{h_{x_i}}(A_g)\delta_{(i,i)}$.

Given $i,j\in \bI_n$ and $g\in\tG$, we define  the additive bijections $\gamma^{(ij)}_g:D^{(i)}_{g\m}M_{ij}\to D^{(i)}_{g} M_{ij}$ by the rule 
\begin{align}\label{gamma-skew}
	\gamma^{(ij)}_g\big(\af_{h_{x_i}}(a_{g\m})\delta_{(j,i)}\big)=\af_{h_{x_i}}(\alpha_g(a_{g\m}))\delta_{(j,i)},\quad a_{g\m}\in A_{g\m}.	
\end{align}
In the rest of this section, we shall denote $\gamma_{g}=\big\{\gamma^{(ij)}_g:D^{(i)}_{g\m}M_{ij}\to D^{(i)}_{g} M_{ij}\big\}_{i,j\in \bI_n}.$ 

\begin{lem1}\label{satisfies-cond} With the notations above $\mcD=\big(\{\alpha^{(i)}\}_{i\in \bI_n}, \{\gamma_{g}\}_{g\in\tG}\big)$ is a datum for $R.$
\end{lem1}
\begin{proof} Let $g\in \tG$ and $\mcI_g=\big\{D^{(i)}_g\,:\,i\in\bI_n\big\}$. We shall first show that  $R$ is $\mcI_g$-symmetric. Indeed, 
	let $i,j\in \bI_n$ and $g\in \tG$. Notice that  
	\begin{align*}
		D^{(i)}_{g}M_{ij}	&=\af_{h_{x_i}}(A_g)\delta_{(i,i)} M_{ij}= \af_{h_{x_i}}(A_g)\delta_{(i,i)} A_{x_i}\delta_{(j,i)}\\
		&\overset{\mathclap{\eqref{eq:structure-module-skew}}}{=}\af_{h_{x_i}}(A_g)\delta_{(j,i)}\overset{\mathclap{\eqref{eq:structure-module-skew}}}{=} M_{ij}\af_{h_{x_j}}(A_g)\delta_{(j,j)}\\
		&=M_{ij}D^{(j)}_{g},
	\end{align*}
	which implies that \eqref{condi-ideals} holds. 
	On the other hand, it is easy to check that $\gamma_g^{(ij)}$ satisfies \eqref{cond0}, for all $i,j\in \bI_n$ and $g\in \tG$. It is clear that $\gamma^{(ii)}_g=\alpha^{(i)}_g$. Since $\alpha$ is group-type, we have that $\af_{h_{x_i}}(A_x)=A_{x_i}$. Thus, using that $\alpha_e=\id_{A}$, we obtain that $\gamma^{(ij)}_e=\id_{M_{ij}}$ and hence \eqref{cond1} follows. In order to prove \eqref{cond2}, consider $u=\alpha_{h_{x_i}}(a_{g\m})\delta_{(k,i)}\in \alpha_{h_{x_i}}(A_{g\m})\delta_{(k,i)}$ and $v=\alpha_{h_{x_k}}(b_{g\m})\delta_{(j,k)}\in \alpha_{h_{x_k}}(A_{g\m})\delta_{(j,k)}$ and note that
	\begin{align*}
		uv&=\alpha_{h_{x_i}}(a_{g\m})\delta_{(k,i)}\alpha_{h_{x_k}}(b_{g\m})\delta_{(j,k)}\\
		&\overset{\mathclap{\eqref{eq:structure-module-skew}}}{=} \alpha_{h_{x_i}}(a_{g\m}) \alpha_{h_{x_i}}(b_{g\m})\delta_{(j,i)}\\
		&=\alpha_{h_{x_i}}(a_{g\m}b_{g\m})\delta_{(j,i)}.
	\end{align*}
	Consequently	
	\begin{align*}
		\gamma_g^{(ij)}(uv)&= \alpha_{h_{x_i}}(\alpha_{g}(a_{g\m}))\alpha_{h_{x_i}}(\alpha_{g}(b_{g\m}))\delta_{(j,i)}\\
		&\overset{\mathclap{\eqref{eq:structure-module-skew}}}{=} \alpha_{h_{x_i}}(\alpha_{g}(a_{g\m}))\delta_{(x_k,x_i)}\alpha_{h_{x_k}}(\alpha_{g}(b_{g\m}))\delta_{(j,k)}\\
		&=\gamma_g^{(ik)}(u) \gamma_g^{(kj)}(v),
	\end{align*}
	which implies \eqref{cond2}. Observe that \eqref{cond03} holds because
	\begin{align*}
		\gamma_{h\m}^{(ij)}\left(\alpha_{h_{x_i}}(A_{g\m})\delta_{(j,i)} \cap \alpha_{h_{x_i}}(A_{h})\delta_{(j,i)}\right)&=\gamma_{h\m}^{(ij)}\left(\alpha_{h_{x_i}}(A_{g\m}\cap A_h)\delta_{(j,i)}\right)\\
		&=\alpha_{h_{x_i}}\left(\af_{h\m}\big(A_{g\m}\cap A_h\big)\right)\delta_{(j,i)}\\
		&\subset\alpha_{h_{x_i}}\left(A_{(gh)\m}\right)\delta_{(j,i)},
	\end{align*}
	for all $g,h\in \tG$, $i,j\in \bI_n$. Finally, let $b=\alpha_{h_{x_i}}(a)\delta_{(j,i)}$, where  $a\in\af_{h\m}(A_{g\m}\cap A_h)$. Then
	\begin{align*}
		\gamma_{g}^{(ij)}\big(\gamma_{h}^{(ij)}(b)\big)&=\alpha_{h_{x_i}}(\af_g(\af_h(a)))\delta_{(j,i)}\\
		&=\alpha_{h_{x_i}}(\af_{gh}(a))\delta_{(j,i)}=\gamma_{gh}^{(ij)}\big(b\big),
	\end{align*}
	and hence \eqref{cond4} is valid.
\end{proof}

\begin{prop1} \label{prop-aux}
	Consider the family of pairs $\gamma=\left(I_g,\gamma_g\right)_{g\in \tG}$, where  
	\[I_g=\big(\af_{h_{x_i}}(A_g)\delta_{(j,i)}\big),\quad \gamma_g=\big(\gamma^{(ij)}_{g}\big),\]
	and $\gamma^{(ij)}_{g}$ is given by \eqref{gamma-skew}, for all $i,j\in \bI_n$ and $g\in \tG$. Then $\gamma$ is a partial action of $\tG$ on the generalized matrix ring $R$.
\end{prop1}
\begin{proof}
	It follows from Lemma \ref{satisfies-cond} and Theorem \ref{teo-par}.	
\end{proof}
Now we shall prove the main result of this subsection. To that end, recall that  the family of pairs $\varepsilon=(C_g, \varepsilon_g)_{g\in \tG}$ given by \eqref{pa-oa} is a partial action of the isotropy group $\tG=\G(x)$ on $A*_{\beta}\G_0^2$.

\begin{teo1} \label{teo-skew-iso}
	The ring isomorphisms
	$ A*_{\af}\G\simeq (A*_{\beta}\G_0^2)*_{\varepsilon} \tG\simeq  R*_{\gamma} \tG$
	hold.
\end{teo1}
\begin{proof}
	The left side isomorphism follows from  \cite[Thm 4.4]{BPP} and it is given explicitly in \eqref{teo-skew}. For the other one, note that  $\psi:A*_{\beta}\G_0^2 \to B*_{\theta}\bI_n^2$ given in \eqref{lem-iso-skews} and $\theta:B*_{\theta}\bI_n^2\to R$ given in \eqref{prop-gmr} are ring isomorphisms. Then $\varphi=\theta\circ\psi$ is a ring isomorphism from $A*_{\beta}\G_0^2$ to $R$ that satisfies $\varphi(C_g)=I_g$ and $\big(\varphi|_{C_g}\big)\circ\varepsilon_g=\gamma_{g}\circ \big(\varphi|_{C_{g\m}}\big)$, for all $g\in \tG$. Thus $\varepsilon$ and $\gamma$ are equivalent  partial actions  and consequently $(A*_{\beta}\G_0^2)*_{\varepsilon}\tG\simeq R*_{\gamma}\tG$.
\end{proof}

\begin{remark}{\rm By the proof of Theorem \ref{teo-skew-iso}, the partial group actions $\varepsilon$  and $\gamma$ are equivalent, where $\varepsilon$ is given in \cite[ Lemma 4.3]{BPP} and $\gamma$ is defined in Proposition \ref{prop-aux}. Hence, Lemmas 4.1, 4.2 and 4.3 of \cite{BPP} are  particular cases of the general construction presented in the previous subsection. }
\end{remark}

\subsection{An example on $2$-by-$2$ generalized matrix ring} Let $B$ be a ring and $K_1,\,K_2$ be ideals of $B$. Suppose that $K_i=Be_{i}$, where $e_i$ are central idempotents of $B$, $i=1,2$. We consider the generalized matrix ring 
$$ R = \left( \begin{matrix} B & K_1\\ K_1 & B \end{matrix} \right).$$
In this case, $M_{11}=R_1=B=R_2=M_{22}$ and $M_{12}=K_1=M_{21}$.
Let $\theta=\{\theta_g\}_{g\in \tG}$ be a global action of a group $\tG$ on $B$ such that $K_1$ is $\theta$-invariant, that is, $\theta_{g}(K_1)\subset K_1$, for all $g\in\tG$. In this case, $\theta_{g}(K_1)=K_1$ and we have that $\theta_g|_{K_1}:K_1\to K_1$ is a ring isomorphism. In particular, $\theta_g(e_1)=e_1$, for all $g\in \tG$. \smallbreak

Now, we consider $\beta^{(i)}=\theta$, which is a global action of $\tG$ on $R_i=B$, for  $i=1,2$. Also, we establish that
\[{\beta}_{g}=\big\{{\beta}^{(ij)}_g:M_{ij}\to M_{ij}\big\}_{i,j\in \bI_2},\text{ where } {\beta}^{(ii)}_g=\beta^{(i)}_g=\theta_g,\,\,{\beta}^{(ij)}_g=\theta_g|_{K_1},\] for all $g\in \tG$ and $i,j\in\bI_2$, $i\neq j$. It is clear that  the elements of ${\beta}_{g}$ satisfy the conditions \eqref{glo-con1} and \eqref{glo-con2}. Hence
$\tilde{\mcD}=\Big(\{\beta^{(i)}\}_{i\in \bI_n}, \{{\beta}_{g}\}_{g\in\tG}\Big)$ is a global datum for $R$. From $\S\,$\ref{subsec-4-3}, we obtain a global action $\beta$ of $\tG$ on $R$ given explicitly by $\beta_g\big((r_{ij})\big)=\big(\theta_g(r_{ij})\big)$, for all $(r_{ij})\in R$ and $g\in \tG$. \smallbreak

Consider $\mathcal{J}=\{J_1,J_2\}$ where  $J_1=J_2=K_2=Be_2$. Since $K_2K_1=K_1K_2=Be_1e_2$ and $K_2B=BK_2=K_2$, it follows that $R$ is $\mathcal{J}$-symmetric. Thus, 
$$ J = \left( \begin{matrix} K_2 & K_1K_2\\ K_1K_2 & K_2 \end{matrix} \right)=\left( \begin{matrix} Be_2 & Be_1e_2\\ Be_1e_2 & Be_2 \end{matrix} \right)$$
is an ideal of $R$. Hence, we have a partial action  $\alpha=(\alpha_g,J_g)_{g\in \tG}$ of $\tG$ on the generalized matrix ring $J$ given by $J_{g}=J\cap \beta_g(J)$ and $\alpha_g=\beta_{g}|_{J_{g\m}}$, for all $g\in \tG$. Since $\theta_g(e_1)=e_1$, it follows that 
$$J_g= J \cap \beta_g(J)=\left( \begin{matrix} Be_2\theta_g(e_2) & Be_2\theta_g(e_2)e_1\\ Be_2\theta_g(e_2)e_1 & Be_2\theta_g(e_2) \end{matrix} \right),$$

and $\alpha_{g}:J_{g\m}\to J_g$ is given by the following association: for all $b_i\in B$, $i\in \bI_4$,
\[\left( \begin{matrix} b_1e_2\theta_{g\m}(e_2) & b_2e_2\theta_{g\m}(e_2)e_1\\ b_3e_2\theta_{g\m}(e_2)e_1 & b_4e_2\theta_{g\m}(e_2) \end{matrix} \right)\mapsto
\left( \begin{matrix} \theta_g(b_1)e_2\theta_g(e_2) & \theta_g(b_2)e_2\theta_g(e_2)e_1\\ \theta_g(b_3)e_2\theta_g(e_2)e_1 & \theta_g(b_4)e_2\theta_g(e_2) \end{matrix} \right).\]
Particularly,  $\alpha^{(i)}=\big(\alpha^{(i)}_g,D^{(i)}_g\big)_{g\in \tG}$ given by
\[D^{(i)}_g=Be_2\theta_g(e_2), \quad\alpha^{(i)}_g\big(be_2\theta_{g\m}(e_2)\big)=\theta_g(b)e_2\theta_g(e_2),\qquad b\in B,\,\,i\in\bI_2,\]
is a  partial action of $\tG$ on $J_i=J_{ii}=Be_2$.
Notice that 
\[D^{(i)}_g=e_2\theta_g(e_2)B, \quad D^{(i)}_gM_{11}=D^{(i)}_g= D^{(i)}_gM_{22}, \quad D^{(i)}_gM_{12}= e_2\theta_g(e_2)e_1B=D^{(i)}_gM_{21}, \] 
for all $g\in \tG$ and $i=1,2$. Also, we have the following additive bijections: $\gamma^{(ii)}_g=\alpha^{(i)}_{g}$ and $\gamma^{(21)}_g=\gamma^{(12)}_g: D^{(i)}_{g\m}M_{12}\to D^{(i)}_gM_{12} $ given by
\begin{align*}
	\gamma^{(12)}_g\big(e_2\theta_{g\m}(e_2)e_1b\big)=e_2\theta_{g}(e_2)e_1\theta_g(b),\quad \text{for all } \,b\in B.
\end{align*}
By Proposition \ref{coincide}, $\mcD=\Big(\{\alpha^{(i)}\}_{i\in \bI_n}, \{\gamma_{g}\}_{g\in\tG}\Big)$ is a datum for $J$ and the partial action $\gamma$ of $\tG$ on $J$ associated to $\mcD$ coincides with $\alpha$. 

\subsection{An explicit example} Let $k$ be a unital commutative ring and consider the ring $B=k^n=k\times\ldots\times k$, with $n\in \N$, $n>1$. Let $r$ be an integer number such that $1\leq r< n$. Let $\tilde{e}_j=(0,\ldots,0,1,0,\ldots,0)$ be the idempotent of $B:=k^n$ which has $1$ in the j-th coordinate and zero in the other ones. Consider
\[e_1=\tilde{e}_1+\ldots+\tilde{e}_r,\quad e_2=\tilde{e}_{r}+\ldots+\tilde{e}_{n-1},\quad K_i=Be_i,\quad i\in\bI_2.  \]
Then $K_1\simeq k^r$ and $K_2\simeq k^{n-r}$  are ideals of $B$. We will use the construction above to give an explicit example of a partial action of the cyclic group $\tG=C_{n-r}=\langle g \rangle$ on the generalized matrix
$$ R = \left( \begin{matrix} B &Be_1\\ Be_1 & B \end{matrix} \right),$$
We define the following ring isomorphism on $k^n$:
\[\theta_g\big((a_1,\ldots,a_r,a_{r+1},\ldots,a_n)\big)=(a_1,\ldots,a_r,a_{n},a_{r+1}\ldots,a_{n-1}),\quad a_j\in k,\,j\in \bI_n.\]

Observe that $\theta_g(e_1)=e_1$ and $\theta_g(e_2)=\tilde{e}_{r}+\tilde{e}_{r+2}+\ldots +\tilde{e}_{n}$ and $e_1e_2=\tilde{e}_r$. Moreover, $\theta=\{\theta_{g^i}\}_{i=1}^{n-r}$, given by $\theta_{g^i}=(\theta_{g})^{i}$, is a global action of $C_{n-r}$ on $B$. By the previous subsection, we have a partial action $\gamma=\left(J_g,\gamma_g\right)_{g\in \tG}$ of $C_{n-r}$ on the generalized matrix ring
\[J =\left( \begin{matrix} Be_2 & B\tilde{e}_r\\ B\tilde{e}_r & Be_2 \end{matrix} \right).\]
Notice that $\theta_{g^i}(e_2)=\tilde{e}_{r}+\ldots+\tilde{e}_{r+i-1}+\tilde{e}_{r+i+1}+\ldots+\tilde{e}_{n}$, for all $i\in\bI_{n-r-1}$. We will denote  $1_{g^{i}}=e_2\theta_{g^i}(e_2)$. Then
\begin{align*}
	1_{g^{i}}=\tilde{e}_{r}+&\ldots+\tilde{e}_{r+i-1}+\tilde{e}_{r+i+1}+\ldots+\tilde{e}_{n-1},\quad i\in\bI_{n-r-2},\\[.2em]
	&1_{g^{n-r-1}}=\tilde{e}_{r}+\ldots+\tilde{e}_{n-2}.
\end{align*}
Thus,
$$J_{g^i}=\left( \begin{matrix}B1_{g^i} & B\tilde{e}_r\\ B\tilde{e}_r &B1_{g^i} \end{matrix} \right), \quad\text{for all }\, i\in \bI_{n-r-1}.$$
For each $i\in\bI_{n-r-2}$, the ring isomorphism $\gamma_{g^i}:J_{g^{n-r-i}}\to J_{g^i}$ is given by
\[\gamma_{g^i}\left(\left( \begin{matrix} x_1 & y_1\tilde{e}_{r}\\ y_2\tilde{e}_{r} & x_2 \end{matrix} \right)\right)=\left( \begin{matrix} a_1 & y_1\tilde{e}_{r}\\ y_2\tilde{e}_{r} & a_2 \end{matrix} \right),\]
where
\begin{align*}
	&x_t=b^t_0\tilde{e}_{r}+b^t_1\tilde{e}_{r+1}+\ldots+b^t_{n-r-i-1}\tilde{e}_{n-i-1}+b^t_{n-r-i+1}\tilde{e}_{n-i+1}+\ldots+b^t_{n-r-1}\tilde{e}_{n-1}, \\
	&a_t=b^t_0\tilde{e}_{r}+b^t_{n-r-i+1}\tilde{e}_{r+1}+\ldots+b^t_{n-r-1}\tilde{e}_{r+i-1} + b^t_1\tilde{e}_{r+i+1}+\ldots+b^t_{n-r-i-1}\tilde{e}_{n-1},
\end{align*}
and the elements $y_1,\,y_2,\,b^t_j\in k$ with $t=1,2$. Finally, $\gamma_{g}:J_{g^{n-r-1}}\to J_{g}$ is the ring isomorphism given by
\[\gamma_{g}\left(\left( \begin{matrix} z_1 & w_1\tilde{e}_{r}\\ w_2\tilde{e}_{r} & z_2 \end{matrix} \right)\right)=\left( \begin{matrix} b_1 & w_1\tilde{e}_{r}\\ w_2\tilde{e}_{r} & b_2 \end{matrix} \right),\]
where
\begin{align*}
	&z_t=c^t_0\tilde{e}_{r}+c^t_1\tilde{e}_{r+1}+\ldots+c^t_{n-r-2}\tilde{e}_{n-2}, \\
	&b_t=c^t_0\tilde{e}_{r}+c^t_{n-r-1}\tilde{e}_{r+2}+\ldots+c^t_{n-r-2}\tilde{e}_{n-1},
\end{align*}
and the elements $w_1,\,w_2,\,c^t_j\in k$ with $t=1,2$.


\begin{thebibliography}{AAH3}
\bibitem{ADES} F. Abadie, M. Dokuchaev, R. Exel and J.J . Sim\'on, \emph{ Morita equivalence of partial group actions and globalization}.  Trans. Amer. Math. Soc. \textbf{368} (1), 4957–4992  (2016).  

\bibitem{BFP} D. Bagio, D. Fl\^ores and A. Paques, \emph{Partial actions of ordered groupoids on rings}. J. Algebra Appl. \textbf{9}, 501--517 (2010).

\bibitem{BLP} D. Bagio, J. R. Lazzarin and A. Paques, \emph{Crossed products by twisted partial actions: separability, semisimplicity and Frobenius properties}.  Comm.  Alg.
\textbf{38}, 496--508 (2010).


\bibitem{BP} D. Bagio and A. Paques, \emph{Partial groupoid actions: globalization, Morita theory and Galois theory}. Comm. Algebra \textbf{40}, 3658--3678 (2012).	

\bibitem{BPP} D. Bagio, A. Paques, and H. Pinedo, \emph{ On partial skew  groupoid  rings}. Internat. J. Algebra Comput. \textbf{31} (1), 1--17 (2021).
	
\bibitem{budanov} A. V. Budanov, \emph{Ideals of generalized matrix rings}. Sb. Math, \textbf{202}, 1--8 (2011).

\bibitem{D} M. Dokuchaev, {\it Recent developments around partial actions,}  São Paulo J. Math. Sci. 13 (1), 195–247 (2019).



\bibitem{DE} M. Dokuchaev and R. Exel, \emph{Associativity of crossed products by partial actions,
	enveloping actions and partial representations}.  Trans. Amer. Math. Soc. \textbf{357}, 1931--1952 (2005).

\bibitem{DFP} M. Dokuchaev, M. Ferrero and A. Paques, \emph{Partial actions and Galois theory}.  J. Pure Appl. Algebra \textbf{208}, 77--87  (2007).

\bibitem{Exel} R. Exel, {\it Partial dynamical systems, Fell bundles and applications}. Mathematical Surveys and Monographs 224, AMS (2017).

\bibitem {F}C. Faith, {\it  Algebra: rings, modules and categories}. vol. I, Springer-Verlag, Berlin–Heidelberg–New York, (1973).


\bibitem{FJ} B. Ferreira and  A. Jabeen,  {\it Multiplicative maps on generalized $n$-matrix rings}. arXiv:2205.15728, (2022).






\bibitem{Ja} N. Jacobson, {\it Basic algebra II}. Second Edition, Dover Publications, (2009).




	
\end{thebibliography}
\end{document}